\documentclass[a4paper,11pt]{article}
\usepackage[margin=3cm]{geometry}
\usepackage{amsfonts,amsmath,amssymb,amscd,graphicx}
\usepackage[colorlinks]{hyperref}

\newtheorem{theorem}{Theorem}[section]
\newtheorem{corollary}[theorem]{Corollary}
\newtheorem{lemma}[theorem]{Lemma}

\newtheorem{proposition}[theorem]{Proposition}

\newtheorem{definition}[theorem]{Definition}
\newtheorem{remark}[theorem]{Remark}

\numberwithin{equation}{section}

\newenvironment{preuve}[1][]
{\vskip 2mm  \noindent\emph{\bf Proof#1. }}{$\Box$ \vskip 2mm}

\newcommand{\diam}{\text{Diam}}

\newcommand{\Bb}{\mathbb{B}}

\newcommand{\Nn}{\mathbb{N}}
\newcommand{\R}{\mathbb{R}}

\newcommand{\Z}{\mathbb{Z}}
\newcommand{\C}{\mathbb{C}}

\newcommand{\E}{\mathbb{E}}

\newcommand{\id}{\text{Id}}

\newcommand{\ep}{\epsilon}

\newcommand{\om}{\omega}

\let\epsilon=\varepsilon

\begin{document}

\title{\bf Systoles  and Lagrangians of\\ random complex algebraic hypersurfaces}
%\date
\author{\sc Damien Gayet}
\maketitle
\begin{abstract}
Let $n\geq 1$ be an integer, $\mathcal L \subset \R^n$ be  a compact smooth affine real hypersurface, not necessarily connected. 
We prove that there exists $c>0$ and $d_0\geq 1$, such that for any  $d\geq d_0$, any smooth complex projective hypersurface $Z$ in $\C P^n$ of degree $d$  contains at least 
$ c\dim H_{*}(Z, \R)$ disjoint Lagrangian submanifolds diffeomorphic to $\mathcal L$, where $Z$ is equipped with the restriction of the Fubini-Study symplectic form (Theorem~\ref{theolag}). If moreover the connected components of $\mathcal L$ have non vanishing Euler characteristic, which implies that $n$ is odd, the latter Lagrangian submanifolds form an independent family of $H_{n-1}(Z, \R)$ (Corollary~\ref{corolag}). We use a probabilistic argument for the proof (Theorem~\ref{theorem2}) inspired by a result by J.-Y. Welschinger and the author on random real algebraic geometry~\cite{GWL}, together with quantitative Moser-type constructions (Theorem~\ref{concrete}). For $n=2$, the method provides a uniform  positive lower bound for the probability that a projective complex curve in $\C P^2$ of given degree  equipped with the restriction of the ambient metric has a systole of small size (Theorem~\ref{theorem1}), which is an analog to a similar bound for hyperbolic curves given by M. Mirzakhani~\cite{mirzakhani13}. Our results hold in the more general setting of vanishing loci of holomorphic sections of vector bundles of rank between 1 and $n$ tensored
by a large power of an ample line bundle over a projective complex $n$-manifold (Theorem~\ref{coroX}).
\end{abstract}

Keywords: {Systole, complex algebraic curve, complex projective hypersurface, Lagrangian submanifold, random geometry, K\"ahler geometry}

\textsc{Mathematics subject classification  2010}: 60D05, 53D12, 32Q15
%60K35= percolation
% 60G15 = Gaussian processes
% 26E05 = real analytic functions
% 53D12 = Lagrangian submanifolds
% 32Q15 = Kahler geometry
% 	60D05   	Geometric probability and stochastic geometry

%\newpage

\tableofcontents
\section{Introduction}

\subsection{Disjoint Lagrangian submanifolds}
On a compact orientable smooth real surface of genus $g>1$,
there exists $3g-3$ disjoint non-contractible closed curves such that two of them are not isotopic. A natural generalization of this phenomenon in a closed symplectic manifold $(X,\omega)$ is to estimate the possible number of disjoint Lagrangian submanifolds of given diffeomorphism type in $X$. The answer is easy for submanifolds which exist as compact smooth manifolds in in $\R^{2n}$, like the torus, since by the Darboux theorem they can be implemented at any scale in $X$, so there exists an infinite number of them. Moreover, when the submanifold $\mathcal L$ possesses a smooth non-vanishing closed one-form, which is the case for the $n$-torus,  this form produces an infinite number of disjoint Lagrangian graphs in $T^*\mathcal L$, hence by Weinstein tubular neighbourhood there exists an infinite number of disjoint Lagrangian submanifolds close to $\mathcal L$, see Remark~\ref{remarkseveral} below. If this is not the case or if the Euler characteristic of $\mathcal L$ is not zero, then it cannot be displaced by a perturbation as a disjoint submanifold, see \S~\ref{smalllag}.  Furthermore, the classes of a finite family of disjoint such Lagrangian submanifolds with non zero Euler characteristic form an independant family of the ambient homology group of degree half of the dimension of $X$, see Lemma~\ref{totally}.

\paragraph{The main result. }
In this paper, we are interested in smooth
projective complex submanifolds equipped with the restriction of  the ambient Fubini-Study K\"ahler form. They have the same diffeomorphic type, because they can be isotoped through smooth hypersurface. For this latter reason, the Moser trick and the fact that the symplectic form has entire periods shows that they are all also symplectomorphic, see~Proposition~\ref{nuke}.
Moreover, they benefit an interesting homological property: for any degree $d$ hypersurface $Z\subset \C P^n$, 
$$\dim H_* (Z,\R) \underset{d\to \infty}{\sim} \dim H_{n-1} (Z,\R)\underset{d\to \infty}\sim d^n.$$
The first asymptotic is a consequence of the Lefschetz hyperplane theorem~\cite{griffiths} and the second one can be estimated through 
the Euler class of the tangent space of $Z$ through Euler characteristic and  Chern classes. 
The main goal of this paper is to prove the following theorem: 
\begin{theorem}\label{theolag}
Let $n\geq 1$ be an integer and $\mathcal L \subset \R^n$  be a compact smooth real affine hypersurface, not necessarily connected. Then,
there exists $c>0$ such that for any $d$ large enough,   any complex hypersurface $Z$ of degree $d$ in $\C P^n$ contains at least $c\dim H_{*} (Z,\R)$ pairwise disjoint  Lagrangian submanifolds diffeomorphic to  $\mathcal L$. 
\end{theorem} 
In fact, we prove this result in the  more general setting of vanishing loci of holomorphic sections of vector bundles of rank between 1 and $n$ tensored
 by a large power of an ample line bundle over a projective complex $n$-manifold, see Theorem~\ref{coroX}.
\begin{corollary}\label{corolag}
Under the hypotheses of Theorem~\ref{theolag},
\begin{enumerate}
\item 
if for any component $\mathcal L_i$ of $\mathcal L$, $\chi(\mathcal L_i)\not= 0$, then the  classes in $H_{n-1}(Z, \R)$ generated by their Lagrangian copies in $Z$ are linearly independent; 
\item if $\mathcal L$ is simply connected, its Lagrangian copies are not close perturbations of each other. 
\end{enumerate}
\end{corollary}
%Let $\mathcal H_{H}$ the set of 
%$\mathcal L$, such that 
%In other terms, in odd dimensions, if $\mathcal H$ denotes the set of compact smooth affine real hypersufaces with non-vanishing Euler characteristics, 
%there exists 
%when $d$ grows to infinity, an increasing proportion of the dimension of the whole topology of $Z$ 
\begin{remark}
\begin{enumerate}
\item 
Note that  $\chi(\mathcal L)\not=0$ implies that $n$ is odd. 
\item The real projective plane $\R P^2$ is a Lagrangian submanifold of $Z= \C P^2\subset \C P^3$ but cannot be a hypersurface in $\R^3$ since any compact hypersurface of $\R^n$ is orientable. 
\item If $Z=\C P^2\subset \C P^3$, $H_2(Z,\Z)$ is generated
by the class of a complex line 
$[D]$.  The integral of the Fubini-Study K\"ahler form $\omega_{FS}$ over $D$ is positive since $\omega_{FS}$ is positive over complex submanifolds, that is $\langle \omega_{FS},[D]\rangle>0 $. However $\langle \omega_{FS}, [\mathcal L]\rangle = 0$ if $\mathcal L$ is a Lagrangian submanifold, so that 
$H_2(Z, \Z)$ cannot be generated by Lagrangian classes. 
 \end{enumerate}
 \end{remark}
In this paper, we prove  Theorem~\ref{theolag}, which is a special case of  Theorem~\ref{coroX},
through a probabilistic argument, see Theorem~\ref{theorem3}:  if we choose
 at random such a projective hypersurface of given large degree, 
 the probability that the conclusion of the theorem holds is positive. Since the hypersurfaces have the same symplectomorphism type, see Proposition~\ref{nuke}, they all satisfy this property. In a parallel paper~\cite{GD}, we prove this theorem with a deterministic proof based on Donaldson-Auroux method~\cite{donaldson96} and \cite{auroux}. 

\paragraph{Other results on disjoint Lagrangian submanifolds.}
As far as the author of the present work knows, essentially three types of results for disjoint Lagrangian submanifolds have been proved. 
\begin{itemize}
\item  The oldest one concerns Lagrangian spheres that naturally germ from singularities of hypersurfaces by Picard-Lefschetz theory. For instance, S. V. Chmutov~\cite[p. 419]{arnold} proved that there exists a singular projective hypersurface of degree $d$ with $c_n d^n+ o(d^n)$ singular points, with $c_n\sim_{n}\sqrt{\frac{2}{\pi n}}.$ When the polynomial defining this hypersurface is perturbed into a non-singular polynomial, the singularities give birth  to disjoint Lagrangian spheres of the associated smooth hypersurface of the same degree.
\item The second result is due to G. Mikhalkin and uses toric arguments:
\begin{theorem}{\cite[Corollary 3.1]{mikhalkin05}}\label{mikh}
For any $n\geq 2 $ and $d\geq 1$, a $2h^{n-1,0}$-dimensional subspace of $H_{n-1}(Z,\R) $ has a basis represented by embedded
Lagrangian tori and spheres, where $Z$ is any smooth projective hypersurface of $\C P^n$. 
\end{theorem}
Here, $h^{n,0}$ is the geometric genus of $Z$, that is the dimension of space of the global holomorphic $n$-form $H^{n,0}(X)\subset H^n(X,\C)$.  It grows like 
$cd^n$ for some $c>0$, as does the dimension of  $H_{n-1}(Z,\C)$ and $\chi(Z). $ 
For Lagrangian spheres, Theorem~\ref{mikh} is more precise than our Theorem~\ref{theolag}, since 
with our method, for an even dimension $n\geq 3$, we cannot know if our Lagrangian spheres  have non-trivial class in $H_{n-1}$, and since the constant $c$ in our bound should be very small compared to the one given by Mikhalkin, as in  Chmutov's theorem. Moreover, our theorem does not say anything new for tori. On the other hand, Theorem~\ref{theolag} asserts that \emph{any } real affine hypersurface is produced as Lagrangian submanifolds in a large quantity in the projective hypersurfaces of large enough degree, and in odd dimension, with a  simple topological restriction, it generates a uniform proportion of the homology of the complex hypersurface.  Moreover, our result
extends to any projective manifold equipped with any ample line bundle, see Theorem~\ref{coroX}.
\item  The third type of results concerns  upper bounds for the number of disjoint Lagrangian submanifolds (not necessarily spheres), and uses Floer techniques, see for instance~\cite{seidel} for results in open manifolds and a survey for older results of this kind. 
\end{itemize}

\subsection{Random complex projective hypersurfaces}
The smooth projective complex hypersurfaces of a given degree $d$, that is smooth vanishing loci in $\C P^n$ of complex homogeneous degree $d$ polynomials, 
form a very natural family of compact K\"ahler manifolds. On the contrary  to the real projective hypersurfaces, that is the vanishing loci  in $\R P^n$ of real polynomials, for fixed $d$ the complex hypersurfaces have 
the same diffeomorphism type. In particular for $n=2$, the smooth complex hypersurfaces in $\C P^2$  of degree $d$ are compact connected Riemann surfaces of genus 
$$ g_d:= \frac{1}{2}(d-1)(d-2).$$
Moreover, as said before, for any $n$ and $d$, when the complex hypersurfaces are equipped with the restriction of the ambiant K\"ahler form, they all have the same  symplectomorphism type, see Proposition~\ref{nuke}.

In~\cite{SZ}, the authors inaugurated the study of random vanishing loci of complex polynomials in higher dimensions (and zero sets of random holomorphic sections, see section~\ref{general}), studying in particular the statistics of the current of integration over these loci. 
In this paper, we will study the statistics of some metric and symplectic properties of these hypersurfaces equipped with the restriction of the Fubini-Study K\"ahler  metric $g_{FS}$ and form $\omega_{FS}$ on $\C P^n$. More precisely, we will be concerned with systoles and small Lagrangian submanifolds. 
\begin{itemize}
\item {\bf $n=2$.}
A source of inspiration and motivation for this paper in the case where $n=2$ was genuinely probabilistic and provided by M. Mirzhakani's theorem on systoles of random hyperbolic curves~\cite{mirzakhani13}, see Theorem~\ref{mimi}.  One of the two main goals of the present  work is in fact to find an analog of it for random complex projective curves, see Theorem~\ref{theorem1}. 
\item {\bf $n\geq 3$.}
With the methods we use, it happens that in higher dimension the natural generalization of small non-contractible loops are small Lagrangian submanifolds of the random hypersurfaces. Our motivation was nevertheless deterministic. The probabilistic method is partly inspired by the work of J.-Y. Welschinger and the author on random \emph{real} algebraic manifolds~\cite{GWL}, were we proved that any compact affine real hypersurface $\mathcal L$ appears a lot of times as a component of a random large degree real projective hypersurface with a uniform probability, see Theorem~\ref{thGWL}. Note that these components are Lagrangian submanifolds of the complexified hypersurface. In particular,  this implies that any complex hypersurface of large enough degree $d$ contains at least $c\sqrt d^n $ Lagrangian submanifolds diffeomorphic to $\mathcal L$, where $c>0$ does not depend on $d$, see Remark~\ref{GWremark}.  In this paper, we prove an analogous complex and symplectic  result analogous to Theorem~\ref{thGWL}:  that  any compact  real hypersurface appears at least $cd^n$ times  as a small Lagrangian submanifold in a random  complex projective hypersurface with a uniform positive probability, see Theorem~\ref{lagprob}.
We emphasize that this improvement from $\sqrt d^n $ to $d^n$ has an interesting topological implication: when $\chi(\mathcal L)\not=0$, it implies that  these disjoint submanifolds form an independent family of homology classes of a cardinal comparable to the dimension of the whole homology  of the complex hypersurface. As said before, the deterministic Theorem~\ref{theorem1} is a direct consequence of the probabilistic Theorem~\ref{lagprob}.
\item It can be suprising that probabilistic arguments
can have deterministic consequences in this situation. The main explanation is given by Theorem~\ref{theorem2} which shows for any sequence of smaller and smaller balls $B$ of size $1/\sqrt d$ in $\C P^n $, the  Lagrangian of desired diffeomorphism type appears in the intersection of $B$ and the random hypersurface  of degree $d$ with a uniform positive probability. This uniform localization  easily implies the global Theorem~\ref{lagprob}, which says that with uniform probability, a uniform proportion of a packing of $\C P^n$ with disjoint balls of size 
$1/\sqrt d$ contain the wanted Lagrangian. 
It happens that the order $d^n$ of growth of $\dim H_{*} (Z, \R)$ is the same order of the number of these packed small balls. 
This result itself implies immediatly the deterministic consequence. 
\item Finally,
 using the universality of peak sections on K\"ahler manifolds equipped with ample line bundles, or the asymptotic (in the degree $d$) universality of the Bergmann kernel, we will explain that analogous results can be proved in this general setting, see the probabilistic Theorem~\ref{theorem3} and the deterministic Theorem~\ref{coroX}.
\end{itemize}

Let us define the measure  on the space of complex polynomials used in~\cite{SZ} and in this paper.
Let $$H_{d,n+1}:= \C_{hom}^d[Z_0, \cdots ,Z_n]$$ be the space of complex homogeneous polynomial in $n+1$ complex variables.
Its dimension equals $n+d\choose n$. For $P\in H_{d,n+1}$, 
denote by $Z(P) \subset \C P^n$ its projective vanishing locus.  
For $P$ outside a codimension 1 complex subvariety of $H_{d,n+1}$,  $Z(P)$ is a smooth complex hypersurface.
%For $n= 2$, they are complex curves 
% of genus $$g_d:=\frac{1}{2}(d-1)(d-2),$$ see~\cite{griffiths} for instance.
Since for transverse polynomials $P,Q$, $Z(P) = Z(Q)$ is equivalent to  $P= \lambda Q$ for some $\lambda \in \C^*$, the space of degree $d$ hypersurfaces  has the dimension of $H_{d,n+1}$ minus one. For $n=2$ this is $\frac{1}{2}d(d+3)\underset{g_d\to \infty}{\sim} g_d$. Note that for the hyperbolic curves, the complex moduli
space has dimension $3g-3$.  
 There exists a natural 
Hermitian product  on $H_{d,n+1}$ given by
$$ \forall P,Q\in H_{d,n+1}, 
\langle P, Q\rangle = \int_{\C P^n} h_{FS}(P,Q) d\text{vol}_{g_{FS}},$$
where $$h_{FS}(P,Q)([Z]) = \frac{P(Z)\overline{ Q(Z)}}{|Z|^{2d}}$$
and  $g_{FS}$ denotes the Fubini-Study metric on $\C P^n$.
Recall that the latter is the quotient metric induced by the 
projection $\C^{n+1}\supset \mathbb S^{2n+1}\to \C P^n$
and the standard round metric on the sphere.  
Then, 
the monomials 
\begin{equation}\label{monomials}
\left(\sqrt{\frac{(d+n)!}{{i_0!\cdots i_n!}}}Z_0^{i_0}\cdots Z_n^{i_n}\right)_{\sum_{k=0}^n i_k =d}
\end{equation}
form an orthonormal basis of $\C_{hom}^d[Z_0, \cdots ,Z_n]$, see the end of the proof of  Lemma~\ref{norm}.
This Hermitian product induces a Gaussian probability  measure on $H_{d,n+1}$. In other terms, 
we choose 
\begin{equation}\label{poly}
P= \sum_{i_0+\cdots +i_n = d} a_{i_0\cdots i_n} \sqrt{\frac{(d+n)!}{i_0!\cdots i_n!}}Z_0^{i_0} \cdots Z_n^{i_n}
\end{equation}
with i.i.d. Gaussian  coefficients $a_I\in \C$ such that $\Re a_I \sim N(0,1)$
and $\Im a_I\sim N(0,1)$ and are independent. We denote by $\mathbb P_d$ the measure.

\subsection{Systoles of random  projective curves}\label{systolesr}

Let $(X,h)$ be a  compact smooth real manifold  equipped with a metric $h$.
In~\cite{mirzakhani13}, M. Mirzakhani studied probabilistic
aspects of metric parameters of $(X,h)$, when 
$(X,h)$ is taken at random in $\mathcal M_g$, the moduli
space of hyperbolic genus $g$ compact Riemann surfaces. 
This moduli space is equipped with a natural 
symplectic form, the Weil-Petersson form, hence a volume form, for
which $\mathcal M_g$ has a finite volume, and
 provides a natural probality measure $\mathbb P_{WP,g}$ on it, see~\cite{mirzakhani13}.
Denote by 
\begin{itemize}
\item $\ell_{sys}(X)$ the least length
 of the noncontractible loops in $(X,h)$.
 \end{itemize}
M. Mirzakhani proved  the following theorem:
\begin{theorem}\cite[Theorem 4.2]{mirzakhani13}\label{mimi}
There exists $\ep_0>0$ and $0<c<C$ such that 
for any $\ep\leq \ep_0 $ and every $g\geq 2$,
$$  \ c \ep^2 \leq \mathbb P_{WP,g} \big[X\in \mathcal M_g \ |\  \ell_{sys}(X) < \ep
\big]\leq C \ep^2. $$ 
\end{theorem}
%\subsection{Random algebraic complex hypersurfaces}
\begin{figure}
$$\begin{array}{|c|c|c|}
\hline
&&\\
 \text{Parameters} &  & \\
 \text{of the surfaces of genus }g & \text{Hyperbolic surfaces} & \text{Planar algebraic curves} \\&& \\ \hline  \hline 
 &&\\
 \text{Dimension of the moduli space} & \underset{g\to\infty}{\sim}  3g & \underset{d\to\infty}{\sim}  g \\
\text{Curvature}& -1 &\in ]-\infty, 2] \ \cite{ness}\\
\text{Volume} & \underset{g\to\infty}{\sim}  4\pi g &  \underset{g\to\infty}{\sim}  4\pi g\\
\text{Diameter} & \in ]0,+\infty[ & \in [c, C g^{5/2}]\ \cite{feng_schumacher_1999}\\
&&\\
\hline
\end{array}$$
\caption{Deterministic parameters of the two different models of real surfaces, the Weil-Petersson one with hyperbolic surfaces, and the Fubini-Study model with complex algebraic curves equipped with the induced rescaled induced metric $\sqrt {2\pi d} g_{FS}$ on $\C P^2$.}\label{tab1}
\end{figure}
\begin{figure}
$$\begin{array}{|c|c|c|}
\hline
&&\\
 \text{Parameters} & \text{Hyperbolic surfaces} & \text{Planar algebraic curves} \\
 \text{of the surfaces of genus }g   & \text{Weil-Petersson measure} & \text{Fubini-Study measure} \\ &&\\ \hline  \hline 
 &&\\
\text{Curvature}
& -1  &\mathbb E_d \big(K(x)|x\in C\big) \asymp -1 \\
\text{Diameter} & \mathbb P_{WP, g} (\diam \geq 40\log g) \underset{g\to\infty}{\to}  0 \ \cite{mirzakhani13}& ? \\
\text{Systole} & \mathbb P_{WP, g} (\ell_{sys}\leq \ep) \asymp \ep^2 \ \cite{mirzakhani13} & \mathbb P_d (\ell_{sys}\leq \ep) \geq \exp(- \frac{c}{\ep^{6}})\  [\text{this paper]}\\
&&\\
&& \\
\hline
\end{array}$$
\caption{Statistics of some metric parameters. The complex algebraic curves equipped with the induced rescaled induced metric $\sqrt {2\pi d} g_{FS}$ on $\C P^2$.}\label{tab2}
\end{figure}
We now introduce a partial analogous result for random projective curves of given degree, with an homological point of view. 
For any $(X,h)$ as above, $\delta>0$ and $c>1$, 
\begin{itemize}
\item 
denote by $N_{sys}(X,\delta,c)$ the 
maximal cardinal of an independent family of classes in 
$H_1(X,\Z)$ such that any class in the family is
represented by 
%we denote by $\ell_\alpha (X, g) $ the least length of the
a  circle of length between $\delta/c$ and $c\delta$. 
\end{itemize}
 Our first main result concerns the systoles of the random complex curves in $\C P^2$:
\begin{theorem}\label{theorem1}There exists a constant $c\geq 1$ and $d_0\geq 1$ such that for every $0<\ep\leq 1$ and every   $d\geq d_0 $,
$$ \exp(-\frac{c}{\ep^6}) < \mathbb P_d \Big[ P\in H_{d,3}\ |   \ N_{sys}\big(Z(P), \frac{\ep}{\sqrt d}, c\big) \geq d^2\exp\big(-\frac{c}{\ep^6}\big) \Big],$$
where $Z(P)$ is equipped with $g_{FS|Z(P)}$. 
In particular, 
$$ \exp\big(-\frac{c}{\ep^6}\big)  < \mathbb P_d \Big[ P\in H_{d,3}\ |\  \ell_{sys} \big(Z(P)\big)  \leq \frac{\ep}{\sqrt d } \Big].$$
\end{theorem}
Theorem~\ref{theorem1} is a particular case of
the more general Theorem~\ref{curve}, which 
holds for random complex curves in a projective complex manifold.
%\begin{corollary}
%There exists a constant $c>0$ and $d_0>0$, $\ep_0>0$, such that 
%for every $\ep<\ep_0$ and every $d\geq d_0,$
%
%\end{corollary}
\begin{remark}
\begin{enumerate}
\item Since $\dim H_1(Z(P), \R) = 2g_d\sim_{d} d^2$, the first assertion of this theorem proves that with uniform probability, 
there exists a basis of $H_1(Z(P), \R)$ such that a uniform proportion of its members is represented by 
a loop of size less than $\ep/\sqrt d$.
\item 
If we want to compare the Fubini-Study model with the Weil-Petersson model, 
we would like that the volumes equal at given genus. This implies that
the metric in the projective setting has to be rescaled with a $\sqrt d$ factor.
In this case  the size  estimates given by Theorem~\ref{theorem1} become similar to the lower bound of Theorem~\ref{mimi}.
Note however that, although our bound is uniform in $d$ or $g_d$ as in~\cite{mirzakhani13},  the dependence in $\ep$ is very bad compared to Mirzakhani's bound. 
\item In fact, for any $x\in \C P^n$, with the same probability, a noncontractible loop lies in $Z(P)\cap B(x,\ep/\sqrt d)$, see Theorem~\ref{theorem2}.
\end{enumerate}
\end{remark}

\noindent
{\bf Other metric parameters. } For the reader's convenience, we present some known results
for other metric properties of the projective curves. Figure~\ref{tab1}. and~\ref{tab2}. compares deterministic and
probabilistic observables for  the Weil-Petersson and  Fubini-Study models.
\begin{itemize}
\item {\bf Volume. } By the Wirtinger theorem, any curve of degree $d$ in $\C P^2$ (and any degree $d$ hypersurface of $\C P^n$) has a volume equal to $d$, see~\cite{griffiths}. By the Gauss-Bonnet theorem, for any hyperbolic curve of genus $g$, the volume equals $2\pi(2g-2)$. Hence, for $n=2$, for comparison with the Weil-Petersson model, 
we should rescale the metric $g_{FS}$ on $\C P^2$ by $\sqrt{ 2\pi d}$, so that 
$$\text{Vol}_{\sqrt {2\pi d} g_{FS|Z(P)} } \big(Z(P)\big)= 2\pi d^2\underset{d\to\infty}{\sim}  4\pi g_d.$$ 
\item {\bf Curvature. } By a result by L. Ness~\cite[Corollary p. 60]{ness}, the Gaussian curvature $K$ of a  degree $d$ complex  curve in $\C P^2$ 
 equipped with the induced metric $g_{FS}$ belongs to $]-\infty, 2]$. 
Besides, by the Gauss-Bonnet theorem, the average on $Z(P)$ of $K$
equals $$K_{mean} = - 2\pi \frac{2g_d-2} d \underset{g\to\infty}{\sim}  -2\pi d.$$
We can prove moreover that 
$ \forall x\in \C P^2, \ \mathbb E \big(K(x)| P(x)=0\big) \asymp -d.$
\item {\bf Diameter. } 
Since by the maximum principle there are no compact complex curves in $\C^2$, no algebraic complex curve in $\C P^2$ does exist in a  ball, so that 
 \begin{equation}\label{didi}
 \exists c>0,  \ \forall P\in \bigcup_{d\geq 1} H_{d,n+1}, \  \diam \big(Z(P), g_{FS|Z(P)}\big) \geq c.
  \end{equation}
 F. Bogomolov~\cite{bogo} has proved that the intrisic diameter of planar complex curves is not bounded when the degree
 grows to infinity. However S.-T. Feng and G. Schumacher~\cite{feng_schumacher_1999} showed that for a given degree
there exists an upper bound for the diameter given by: 
$$ \forall d\geq 1, \ \forall P\in H_{d,3}, \ \diam \big(Z(P),  g_{FS|Z(P)} \big) \leq 32\pi g_d^{2}+ o(g_d^2).$$
 It should be possible, like in~\cite{mirzakhani13}, to find a better probabilistic estimate for the diameter, and one can wonder if it is also logarithmic in $d$. 
\end{itemize}
%\begin{remark} M. Gromov proved~\cite[Theorem 0.2B]{gromov} that there exists a universal constant $A>0$, such that for any Riemannian compact manifold of dimension $n$, 
%if the Gaussian curvature $K$ satisfies 
%$ K\geq -\kappa^2$, then
%$ \sum_i b_i \leq A^{1+\kappa \diam(M)}$, where $b_i$ denotes the $i-$th Betti number of $M$.
%In the case of degree $d$ complex curves,
%this gives
%$ \diam (Z) \geq C \frac{\log d }{\kappa}.$
%By the estimates for the curvatures given above, 
%$\kappa$ should be of the order $\sqrt d$, which
%gives a trivial bound compared to (\ref{didi}).
%\end{remark}

\subsection{Small Lagrangian submanifolds of random hypersurfaces}\label{smalllag}

Let $(X^{2n},\omega)$ be a smooth symplectic manifold of dimension $2n$. Recall that $\omega$ is a closed non-degenerate two-form.  A Lagrangian submanifold $\mathcal L$ of $X$
is a $n$-dimensional submanifold such that $\omega_{|TL}$ vanishes. For instance,  a real analytic hypersurface in $\R^n$ is a Lagrangian submanifold of its associated complex extension, which is a K\"ahler manifold for the restricton of the standard K\"ahler form in $\C^n$.

\paragraph{Universal real components. }
In~\cite{GWL},  J.-Y. Welschinger and the author of the present paper studied random \emph{real} projective hypersurfaces, that is the real loci of random elements of $\R H_{d,n+1}$, the space of real homogeneous  polynomials in $(n+1)$ variables and of degree $d$. The measure was the complex Fubini-Study (\ref{poly}) restricted to $\R H_{d,n+1}$. In the litterature, this measure is often called the \emph{Kostlan measure}. 
 Let $\mathcal L \subset \R^n$ be any compact smooth real  hypersurface.
 For any  real homogeneous polymial $P$, let $Z_{\R}(P):= Z(P)\cap \R P^{n+1},$ and denote by 
 \begin{itemize}
 \item 
 $ N_\R(\mathcal L, Z_\R(P)) $ be the number of disjoint balls  in $\R P^n$  such that for any such ball $B$, $B\cap Z_{\R}(P)$ contains a  submanifold $\mathcal L'$ diffeomorphic to $\mathcal L$.
 \end{itemize}
  
\begin{theorem}{\cite[Theorem 1.2]{GWL} and \cite[Theorem 2.1.1]{GW}}\label{thGWL}
Let $n\geq 1$ and $\mathcal L\subset \R^n$ be any compact smooth  hypersurface, not necessarily connected. 
Then there exists $c>0$ and $d_0$,  such that for every $d\geq d_0$,
$$
\ 
c<
\mathbb P_d 
\Big[
 P\in \R H_{d,n+1} \ | \  N_\R\big(\mathcal L, Z_{\R}(P)\big) >c\sqrt d^n 
\Big].$$
\end{theorem}
\begin{remark}  \label{GWremark}
\begin{enumerate}
\item\label{GWremark1}
Note that this theorem has a deterministic corollary, using the same argument given in this paper: any compact real affine hypersurface appears at least $c\sqrt d^n$ times as disjoint Lagrangian submanifolds in any complex projective hypersurface of high enough degree. Indeed, the real part of a complex hypersurface defined over the reals is Lagrangian for the restriction for the Fubini-Study K\"ahler form, and the complex projective hypersurface are all symplectomorphic. 
\item In \cite{GDR}, the author constructed real hypersurfaces with  $c\sqrt d^n$ real spheres. The same proof, replacing a polynomial
vanishing on a sphere by another polynomial gives the same corollary as the latter. 
Theorem~\ref{theolag} gives a $cd^n$ lower bound, which is of the order of $\dim H_{*}(Z(P), \R)$ when $d$ grows to infinity.
\item 
In fact,  Theorem~\ref{thGWL} holds in the more general context of 
K\"ahler compact manifolds with holomorphic line bundles equipped with real structures, see~\cite{GWL}.
\end{enumerate}
\end{remark}

\paragraph{Universal Lagrangian submanifolds. }
We turn now to a complex and Lagrangian analog of this theorem. As before, let $\mathcal L\subset \R^n$ be a compact smooth real hypersurface. For any compact symplectic manifold $(Z,\omega,h)$ equipped
with a Riemannian metric $h$, any $\delta>0$ and $c\geq 1$, 
%\item 
%  we denote by $N_{Lag} (\mathcal L, X,\delta) $ the number of pairwise disjoint  open sets of $X$, distant from each other from at least $\delta$ and containing a Lagrangian submanifold diffeomorphic to $ \mathcal L$
%  and of diameter bounded by $\delta$. 
\begin{itemize}
\item denote by $N_{Lag}(\mathcal L, Z, \delta, c)$ the number of pairwise disjoint  open sets containing a Lagrangian submanifold  
$\mathcal L'$ 
diffeomorphic to $\mathcal L$ and satisfying:
\begin{equation}\label{distlag}
\frac{\delta}{c}\leq 
\diam\big(\mathcal L', h_{|\mathcal L'}\big) \leq c\delta.
\end{equation}
\end{itemize} 
For polynomials,  the following theorem is the main probabilistic result of this paper. It is a particular case of Theorem~\ref{theorem3} below:
\begin{theorem}\label{lagprob}
Let $n\geq 2$,  $\mathcal L\subset \R^n$ be any compact smooth hypersurface, not necessarily connected. Then there exists $c\geq 1, D\geq 1, d_0\geq 1$
 such that  for any $0<\ep\leq 1$ and  $d\geq d_0$ 
$$
\ \exp(-\frac{c}{\ep^{D}})<
\mathbb P_d 
\left[
P\in H_{d,n+1} \ | \ N_{Lag }\Big(\mathcal L, Z(P), \frac{\ep}{\sqrt d}, c\Big) > d^n \exp\big(-\frac{c}{\ep^{D}}\big)
\right],$$
where the metric and the symplectic form on $Z(P)$ are the ones induced by the Fubini-Study metric and symplectic form on $\C P^n$.
Moreover, is $\mathcal L$ is real algebraic, that is if there exists $p\in \R[x_1, \cdots, x_n]$ such that $\mathcal L= Z_\R(p)$, then $D$ can be chosen to be $D=2\deg p$. 
%The same statement holds, replacing  $N_{Lag}$by $  N_{Lag }^{rel} \big(\mathcal L, Z(P), \C P^n , \frac{\ep}{\sqrt d}\big)$.
%Lastly, if $\chi(\mathcal L) \not= 0$, then the classes of these Lagrangian copies of $\mathcal L$ form an independent family of $H_{n-1}(Z(P), \R)$. 
\end{theorem}
\begin{remark}\label{remarkseveral}
\begin{enumerate}
\item In fact, Theorems~\ref{theolag}, ~\ref{theorem1} and~\ref{lagprob} have a higher codimension generalization: instead of taking one unique random polynomial, one can choose $1\leq r\leq n$  random independent polynomials $(P_1, \cdots, P_r)$ of the same degree, and look at their common vanishing locus $Z(P_1, \cdots,  P_r):=\cap_{i=1}^r Z(P_i)\in \C P^n$, which is now almost surely of complex codimension $r$. Then, the same conclusions hold with the following changes: for complex curves (Theorem~\ref{theorem1}), we take $n\geq 2$ instead of $n=2$, and choose $r=n-1$. For Lagrangians (Theorem~\ref{lagprob}), we take $\mathcal L\subset \R^{n-r+1}$ instead of $\mathcal L\subset \R^n$.  However, if $r\geq 2$,  $\mathcal L\subset \R^{n-r+1}$ 
must satisfy a further necessary condition: its normal bundle must be trivial. 
These generalizations are direct consequences of Theorem~\ref{curve} for the curves  and Corollary~\ref{corolagrangien} for the higher dimensions.
\item 
By the Weinstein theorem~\cite{weinstein}, a tubular neighborhood of a  closed Lagrangian submanifold $\mathcal L$ is symplectomorphic to a tubular neighborhood of the zero section in $T^*\mathcal L$, so that the local Lagrangian deformations of $\mathcal L$ can be viewed in $T^*\mathcal L$ as graphs of closed  $1$-forms on $\mathcal L$.  If the form is exact, then it has at least two zeros and the associated graph intersects $\mathcal L$. In particular, if $H^1(\mathcal L, \R)= 0$, $\mathcal L$ cannot locally be deformed as a disjoint Lagrangian submanifold. 
On the other hand, if $\mathcal L$ possesses a closed $1-$form which does not vanish, like the torus, then there exists an infinite number of Lagrangian submanifold of diffeomorphic to $\mathcal L$. See~\cite{farber} for topological  conditions on $\mathcal L$ which imply non-existence of such non-vanishing forms. This remark shows that in the case of spheres, the disjoint Lagrangian spheres produced by Theorem~\ref{theolag} are not small deformations of each others. 
\item Theorem~\ref{lagprob} is the consequence of the more
precise Theorem~\ref{theorem2}, which asserts
that for any sequence of balls centered on a fixed point $x$ in $\C P^n$ and of size $1/\sqrt d$, with uniform probability $\mathcal L$ appears as a Lagrangian submanifold of the random vanishing locus. 
\end{enumerate}
\end{remark}

\subsection{Random sections of a holomorphic vector bundle}\label{general}
There is at least two natural generalizations of Theorems~\ref{theorem1} and~\ref{lagprob}: firstly, we can work in the 
setting of ample holomorphic line bundles over compact K\"ahler manifolds introduced by B. Schiffman and S. Zelditch in~\cite{SZ}. Secondly, we can study the statistics of the vanishing locus of \emph{several} random polynomials or sections, as said in Remark~\ref{remarkseveral}. We present the fusion of the two generalization, as in~\cite{GW}. 
Let $n\geq 1$ and $X$ be a compact  complex $n$-dimensional manifold equipped with an ample holomorphic line bundle $L\to X$, that is there exists a Hermitian metric $h_L$ on $L$ with curvature 
$-2i\pi \omega$, such that $\omega $ is K\"ahler.  
We denote by
$g_\omega$ the associated K\"ahler metric. 
Note that by the Kodaira theorem, $X$ can be holomorphically embedded in $\C P^N$ for $N$ large enough. 
Let $1\leq r\leq n$ be an integer and $E\to X$ be a 
holomorphic vector bundle of rank $r$ equipped with a Hermitian metric $h_E$. 
For any \emph{degree} $d\geq 1$, denote by $H^0(X,E\otimes L^{ d})$ the space of holomorphic sections of $E\otimes L^{\otimes d}$. By the Hirzebruch–Riemann–Roch theorem, 
$$\dim H^0(X,E\otimes L^{\otimes d})\underset{d\to \infty}{\sim} rd^n\int_X \frac{\omega^n}{n!}. $$
Let $d\text{vol}$ be any volume form on $X$, and define for any $d\geq 1$
the Hermitian product on $H^0(X, E\otimes L^d) $:
\begin{equation}\label{SP}
\forall s,t \in H^0(X, E\otimes L^d), \ 
\langle s, t\rangle := \int_X h_{E,L^d} (s, t) d\text{vol},
\end{equation}
where $h^d_{E,L}$ is the Hermitian metric on $E\otimes L^{\otimes d}$ associated to $h_E$ and $h_L$.
Then we associate to this Hermitian product the Gaussian probability
measure $d\mathbb P_d$ on $H^0(X, E\otimes L^d)$. 
In other term, for any $d\geq 1$, choosing an orthnormal basis $(S_{i})_{i\in \{1, \cdots N_d\}}$ of $H^0(X,E\otimes L^{d}), $ where $N_d:= \dim H^0(X,E\otimes L^d)$, 
a random section $s\in H^0(X, E\otimes L^{ d})$ writes
$$ s= \sum_{i=1}^{N_d} a_i S_i,$$
where the complex coefficients $(\Re a_i)_i$ and $(\Im a_i)_i$  are i.i.d. and follow the same normal law $N(0,1)$. 
In the sequel 
\begin{itemize}
\item
$Z(s)$ will denote the vanishing locus in $X$ of $s\in H^0(X,E\otimes L^{d})$, 
\item and the tuple $(n, r, X,L,E,h_L,\omega, g_\omega, h_E, d\text{vol}, (\mathbb P_d)_{d\geq 1})$
will be called an \emph{ample probabilistic model}, and
\emph{ample model},  if no probability is involved.   
\end{itemize}
 By Bertini's theorem, almost surely 
$Z(s)$ is a compact smooth codimension $r$ complex submanifold of $X$. \\

\noindent
{\bf Standard example: the Fubini-Study random polynomial mappings.} For $X= \C P^n$, $E= \C P^n \times \C^r$, $h_E$ the standard metric on $\C^r$, $L= \mathcal O(1)$ the hyperplane bundle, $h_L= h_{FS}$ the Fubini-Study metric, then 
$$H^0(\C P^n, E\otimes L^d) = \big(\C^d_{hom}[Z_0,\cdots, Z_n])^r.$$ 
Moreover, 
the monomials given by~(\ref{monomials}) make this identification an isometry.  In other terms, a random polynomial mapping for the standard stuctures is a $r$-uple of independent random polyomials in 
$H_{d,n+1}$ equipped with the Gaussian measure~(\ref{poly}).\\

\noindent
{\bf Random curves.} When $r=n-1$, 
the vanishing locus of the sections of $H^0(X, E\otimes L^d)$ 
is generically a smooth compact complex curve. When $n=2$ and $r=1$, 
 the adjunction formula shows that their
genus equals
$$ g_d = \frac{1}{2}d^2 \int_X \omega^2 - \frac{1}{2}d\int_X c_1(X)\wedge \omega +1,$$
 where $c_1(X)$ denotes the first Chern class of the surface $X$, see~\cite{griffiths}. Theorem~\ref{theorem1} has
 the following natural generalization:
\begin{theorem}\label{curve} Let $n\geq 2$ be an integer.  Then, there exists a universal 
constant $c\geq 1$ such that the following holds. 
Let $(n, n-1, X,L,E,h_L,\omega, g_\omega, h_E, d\text{vol}, (\mathbb P_d)_{d\geq 1})$
 be an ample probabilistic model. 
Then, there exists $d_0\geq 1$ such that
 for every $0<\ep\leq 1$ and every $d\geq d_0, $
$$ \exp(-\frac{c}{\ep^{6}})< \mathbb P_d \Big[ s \in H^0(X,E\otimes L^{d}) \ | \  N_{sys}\big(Z(s), \frac{\ep}{\sqrt d}, c\big)  > d^n \text{Vol}_{g_\omega} (X) \exp(-\frac{c}{\ep^{6}})\Big] .$$
Here, the metric on $Z(s)$ is the restriction of the K\"ahler metric $g_\omega$ associated to $\omega$. 
\end{theorem}
Recall that $N_{sys}$ is defined in \S~\ref{systolesr}.
Note that the volume involved in the Theorem~\ref{curve}
is the one associated to $g_\omega$ and not
with the arbitrary volume form $d\text{vol}$ used
for the definition of the scalar product~(\ref{SP}).

Theorem~\ref{curve} means that for any degree large enough, with uniform probability in $d$, 
there exists a basis  $H_1(Z(s), \R)$ such that 
a uniform proportion of its elements are represented by 
loops of size bounded by $\ep/\sqrt d$. 
\begin{remark}
It is classical~\cite[Corollary 3.6]{hartshorne2013algebraic} that any compact orientable Riemann surface embeds in $\C P^3$. However, 
a degree $d$ curve in $\C P^3$, that is a holomorphic curve whose class
in $H_2(\C P^3, \Z)$ equals $d[D]$, where $D$ is a line, can have different
topologies,  and it is not known which pairs
of genus and degree do exist, see~\cite[IV, 6]{hartshorne2013algebraic}. Finally, 
if $E$ is of rank $2$, our model of sections of $H^0(\C P^3, E\otimes L^{d})$ only provides 
strict families of curves of the whole set of curves. 
\end{remark}

\paragraph{Lagrangian submanifolds. }
We provide now a similar K\"ahler generalization of Theorem~\ref{lagprob}, 
that is for Lagrangian submanifolds.
%\begin{itemize}
%\item 
Let    $\Sigma$ be a complex submanifold in $\mathbb B\subset \C^n$, and $\mathcal L$ be a compact smooth Lagrangian submanifold of $(\Sigma, \omega_{0|T\Sigma})$. For   any symplectic  manifold $(Z,\omega,h)$ equipped with a metric $h$ and $\delta>0, c>1$, 
\begin{itemize}
\item 
 denote by $N(\Sigma, \mathcal L, Z,\delta,c)$ the maximal number of pairwise disjoint  open sets $\Sigma'\subset Z$  such that $\Sigma'$ contains a Lagrangian submanifold $\mathcal L'$ such that 
 $$(\mathcal L', \Sigma')\sim_{diff} (\mathcal L, \Sigma) \text{ and }
  \frac{\delta}{c}\leq \diam_{\mathcal L'}\mathcal L'\leq c\delta.$$
 \end{itemize}
%\item We recall that for any compact smooth real affine submanifold $\mathcal L$ and any symplectic manifold $(X,\omega, h)$ equipped with a metric, $ N_{Lag} (\mathcal L, X,\delta) $ denotes the number of disjoint balls of radius $\delta$ in $X$ such that for any such ball  contains a Lagrangian submanifold diffeomorphic to $ \mathcal L$. The other number $N^{rel}_{Lag}(\mathcal L, Z, X, \delta)$ denotes 
%\end{itemize}
\begin{theorem}\label{theorem3} Let $n\geq 2$, 
$1\leq r\leq n-1$,  $\Sigma\subset \mathbb B\subset \C^{n} $ be a complex algebraic smooth codimension $r$ submanifold,   let  $\mathcal L\subset \Sigma $ be a compact smooth Lagrangian submanifold of $(\Sigma, \omega_{0|T\Sigma}).$
Then there exist
$c, D\geq 1$  such that the following holds. 
Let $(n, r, X,L,E,h_L,\omega, g_\omega, h_E, d\text{vol}, (\mathbb P_d)_{d\geq 1})$
be an ample probabilistic model.
Then, there exists  $d_0\geq 1$ such that 
 for every $0<\ep\leq 1$ and $d\geq d_0$, 
\begin{eqnarray*}
\exp(-\frac{c}{\ep^{D}})&<& 
 \mathbb P_d 
 \Big[ s \in H^0(X,E\otimes L^{\otimes d}) \ | \ 
 N\big(\Sigma, \mathcal L, Z(s), \frac{\ep}{\sqrt d}, c\big) 
  > d^n \text{Vol}_{g_\omega} (X)\exp(-\frac{c}{\ep^{D}})\Big].
  \end{eqnarray*} 
\end{theorem}
The following corollary proves that any compact smooth real affine codimension $(n-r)$ submanifold with trivial normal bundle appears a large number of times in the random complex codimension $r$ submanifold, with a uniform probability: 
\begin{corollary}\label{corolagrangien} Let $n\geq 2$, 
$1\leq r\leq n-1$,   and  $\mathcal L\subset \R^{n}$ be a compact smooth codimension $r$ submanifold with trivial normal bundle.
Then there exist
$c, D\geq 1$  such that the following holds. 
Let $(n, r, X,L,E,h_L,\omega, g_\omega, h_E, d\text{vol}, (\mathbb P_d)_{d\geq 1})$
be an ample probabilistic model.
Then, there exists  $d_0\geq 1$ such that 
 for every $0<\ep\leq 1$ and $d\geq d_0$, 
\begin{eqnarray*}
\exp(-\frac{c}{\ep^{D}})&<& 
 \mathbb P_d 
 \Big[ s \in H^0(X,E\otimes L^{\otimes d}) \ | \ 
 N_{Lag}\big(\mathcal L, Z(s), \frac{\ep}{\sqrt d}, c\big) 
  > d^n \text{Vol}_{g_\omega} (X)\exp(-\frac{c}{\ep^{D}})\Big].
  \end{eqnarray*} 
    If $\mathcal L$ is  algebraic, one can choose $D$ as the double of  the degree of $\mathcal L$.
\end{corollary}
Recall that $N_{Lag}$ is defined in \S~\ref{smalllag}.
Note that when $r=1$, that is if $\mathcal L$ is a hypersurface,
the condition on its normal bundle is always satisfied. 
Corollary~\ref{corolagrangien}  implies the following generalization of the deterministic Theorem~\ref{theolag}:
\begin{theorem}\label{coroX}
Let $n\geq 2$, $1\leq r \leq n$ and 
$\mathcal L\subset \R^n$  be a compact smooth $(n-r)$-submanifold with trivial normal bundle. Then, there exists $c>0$ such that for any ample model $(n, r, X,L,E,h_L,\omega, g_\omega,  h_E)$, for $d$ large enough, the zero locus of any section $s\in H^0(X,E\otimes L^{d}) $ vanishing transversally contains at least $cd^n\text{Vol}_{g_\om}(X)$ disjoint  Lagrangian submanifolds diffeomorphic to  $\mathcal L$. 
\end{theorem}
Again, by the Lefschetz theorem and a computation with Chern classes,  there exists $c>0$ such that 
$$\forall d\gg 1, \forall s\in H^0(X, E\otimes L^d), \ \dim H_{*}(Z(s),\R)\underset{d\to \infty}{\sim}\dim H_{n-r}(Z(s),\R)\underset{d\to \infty}{\sim} c d^n,$$
see~\cite[Corollary 3.5.2]{GW} for a proof of it with an explicit constant $c$. 
\begin{corollary}\label{corolagX}
Under the hypotheses of Theorem~\ref{coroX},
\begin{enumerate}
\item 
if for every connected component $\mathcal L_i$ of $\mathcal L$, $\chi(\mathcal L_i)\not= 0$, then the  classes in $H_{n-r}(Z(s), \R)$ generated by these disjoint submanifolds are linearly independent; 
\item if the $\mathcal L_i's$ are simply connected, no Lagrangian copy of them can be isotoped to another one as disjoint Lagrangian submanifolds. 
\end{enumerate}
\end{corollary}
% % % % % % % % % % %
\subsection{Prescribed topology in a small ball}\label{smallb}
Theorem~\ref{theorem3} is a consequence of the more precise Theorem~\ref{theorem2} below. 
This theorem is partly inspired by the work of J.-Y. Welschinger and the author. 
For this reason, we recall it. 
In~\cite{GW}, it was proved the following:
\begin{theorem}{\cite[Proposition 2.4.2]{GW}}\label{realpresence}
Let $n\geq 2$ and $1\leq r\leq n$. Then, for any real compact smooth $(n-r)$-submanifold $\mathcal L \subset \R^n$ with trivial normal bundle, 
for any $d$ large enough and any $x\in \R P^n$, with uniform positive probability in $d$, 
the zero set $Z(P)$ of 
a random real polynomial $P\in \R_{hom}^d [X_0, \cdots, X_n]$ 
intersects $B(x,1/\sqrt d)$ along some components, ones of which are diffeomorphic to $\mathcal L$. 
\end{theorem}
This theorem was in fact proved in the general setting of random sections of holomorphic real vector bundles over a projective manifold, see~\cite{GW}.
We begin with an analogous version of Proposition~\ref{realpresence} for smooth complex algebraic affine hypersurfaces $\Sigma\subset \C^n$
containing a Lagrangian submanifold $\mathcal L$. Note that
the latter condition is not a constraint since every symplectic manifold
contains a Lagrangian torus of any small size enough near every point. 
Note that  contrary to the real case, an affine algebraic complex hypersurface is never compact,  and is connected if and only  it is the vanishing locus of an irreducible polynomial. 
Let $n\geq 2$ and $1\leq r\leq n$ be integers,   $\Sigma$ be a complex submanifold in $\mathbb B\subset \C^n$,  $\mathcal L$ be a compact smooth Lagrangian submanifold of $(\Sigma, \omega_{0|T\Sigma})$, and $(n, r, X,L,E,h_L,\omega, g_\omega, h_E, d\text{vol}, (\mathbb P_d)_{d\geq 1})$ be an ample probabilistic model, see \S~\ref{general} for the definition. For any $x\in X$, $\delta>0$ and $C>1$, define:
\begin{itemize}
\item for any $s\in H^0(X, E\otimes L^d)$, $A\big(\Sigma, \mathcal L, Z(s),x,  \delta, C\big)$ denotes the event that there exists a smooth topological ball $B\subset X$ containing $x$ and a Lagrangian submanifold $\mathcal L'$ of $\big(Z(s)\cap B, \omega_{|Z(s)}\big)$, such that   
$$(\mathcal L', Z(s)\cap B)\sim_{diff} (\mathcal L, \Sigma) \text{ and }
  \frac{\delta}{c}\leq \diam (\mathcal L')\leq c\delta.$$
\end{itemize}
Here, the diameter is computed with respect to the induced metric on $\mathcal L'$. The main theorem of this paper is the following:
\begin{theorem}\label{theorem2} Let $n\geq 2$, $1\leq r\leq n-1$ be integers, 
$\Sigma\subset \overline{\mathbb B} \subset \C^n$ be a smooth complex algebraic $(n-r)$-submanifold, and 
$\mathcal L\subset \Sigma$ be a compact smooth Lagrangian submanifold of $(\Sigma, \omega_{0|T\Sigma})$. Then there exists $c\geq 1$, such that for any ample probabilistic model $(n, r, X,L,E,h_L,\omega, g_\omega, h_E, d\text{vol}, (\mathbb P_d)_{d\geq 1})$, there exists $d_0\geq 1$ such that for any $0<\ep\leq 1$ and for any $x\in X,$
$$ \forall d\geq d_0, \ 
\exp(-\frac{c}{\ep^{D}})\leq 
\mathbb P_d 
\Big[
   A\big(\Sigma, \mathcal L, Z(s), x,\frac{\ep}{\sqrt d }, c\big)
\Big].
$$
\end{theorem}
This theorem implies quickly Theorem~\ref{theorem3}, see below. 
In fact, the same result holds for affine real hypersurfaces, not only Lagrangians, as in Corollary~\ref{corolagrangien} and Theorem~\ref{realpresence}.

\subsection{Ideas of the proof of the main theorems}  
We present the strategy of the proofs of Theorems~\ref{theorem2}, ~\ref{curve}  and~\ref{theorem3} for $r=1$, $\epsilon=1$ and for polynomials. The proof holds on two main tools, namely 
the barrier method for proving uniform probability of some local topological event, together with a quantitative Moser-type construction to make this event symplectic and Lagrangian.
The barrier method
was used for instance in a real deterministic context
in~\cite{GDR} to construct a lot of small spheres in the real part of 
holomorphic or symplectic Donaldson hypersurfaces. 
In probabilistic contexts similar to this present work, it was used for instance  in~\cite{NS1}
to produce small components of the vanishing locus of a  random function  with uniform probability, and
in \cite{GWL} to produce small components with prescribed diffeomorphism types. 
The proof of the main Theorem~\ref{theorem2} is roughly the following:
\begin{itemize}
\item 
fix a point $x\in \C P^n$ and
 choose for any $d$ 
a polynomial $Q_{x,d}$ vanishing 
along a hypersurface  $Z(Q_{x,d})$ intersecting  $B(x,1/\sqrt d)$
along a hypersurface 
diffeomorphic to $\Sigma$. Here, $1/\sqrt d$ of the natural scale for Fubini-Study or Kostlan measures.
The easiest way to do this is to rescale for every $d$ the same polynomial in an affine chart centered on $x$.
\item Then,  for small enough perturbations, the perturbed polynomial  still
vanishes in $B(x,1/\sqrt d)$  along a hypersurface
isotopic to $\Sigma$. If the allowed pertubation can be quantified, typically
when the two-point correlation function of the random function 
converges locally to a universal random function after rescaling, 
one can prove that  with a uniform positive
probability, a random polynomial of degree $d$ vanishes  
in the sequence of balls $B(x,1/\sqrt d)$ along a hypersurface diffeomorphic to $\Sigma$. 
In our case, we specialize this method in two different ways, depending on the dimension $n$ of the ambient space.
\item For $n= 2$ (Theorem~\ref{curve}), 
we choose 
$\Sigma\subset \mathbb B \subset \C^2$ to be a complex curve of degree 3, 
hence a torus without three small disks. Then a circle whose class in $H_1(\Sigma, \R)$ 
is non-trival will still be non-trivial in $H_1(Z(P), \R)$.
\item For $n\geq 3$ (Theorem~\ref{theorem3}), 
In normal affine complex coordinates on the small ball $B(x,1/\sqrt d)$, the Fubini-Study form equals the standard form at  $x$, so that 
the local implementation in $\C P^n$  of $\mathcal L$  is almost Lagrangian in $\big(Z(Q_{x,d}), \omega_{FS|TZ(Q_{x,d})}\big)$. Since the perturbation of $Q_{x,d}$ by a random polynomial is complex and not real,  there is no natural way to follow $\mathcal L$ as a Lagrangian perturbation in the perturbed vanishing locus $\Sigma'$. The classical way to deform  object of symplectic nature, like the Lagrangians, is the Moser method. We reprove it in our particular situation, but with a  quantitative point of view (Theorem~\ref{concrete}). Thanks to the latter the method keeps the perturbation of $\mathcal L$ inside the small ball, so that this small Lagrangian displacement happens with uniform probability. 
These points provide the idea of the proof of Theorem~\ref{theorem2}.
Note that the quantiative Moser trick is needed only for dimensions $n\geq 3$ and not for our result on systoles. 
\item Then, Theorems~\ref{curve} and~\ref{theorem3} are direct consequences of Theorem~\ref{theorem2}: choosing in $\C P^n$ a maximal set of  small disjoint balls, automatically with uniform probability, 
 at least $cd^n$ of these balls intersect
$Z(P)$ along  a component diffeomorphic to $\Sigma$ and contain a Lagrangian copy of $\mathcal L$ with the good diameter. 
\end{itemize}

\paragraph{\bf Organization of the paper. } In section~\ref{poc}, we  we assume Theorem~\ref{theorem2} and  we give the proofs of its consequences presented above. In section~\ref{quant},
we give a quantitative version of the Moser trick. This part is deterministic. In section~\ref{pot}, 
we give the proof of Theorem~\ref{theorem2}.

\paragraph{\bf Aknowledgments.} The author would like to thank Denis Auroux,
Vincent Beffara, Sylvain Courte, Laura Monk, Alejandro Rivera, and Jean-Yves Welschinger for valuable discussions.  The research leading to these results has received funding from the French Agence nationale de la recherche, ANR-15CE40-0007-01 and from the European Research Council project
ALKAGE, contract 670846 from Sept.~2015. 

\section{Direct proofs}\label{poc}
In this section we assume Theorem~\ref{theorem2} and  we give the proofs of its consequences.
\subsection{From local to global }

\begin{preuve}[ of Theorem~\ref{theorem3}]
We follow the proof given in~\cite[\S 2.5]{GW}. 
Let $c\geq 1$ be  given by Theorem~\ref{theorem2},
and let $(n, r, X,L,E,h_L,\omega, g_\omega, h_E, d\text{vol}, (\mathbb P_d)_{d\geq 1})$ be an ample probabilistic model.
Let $\Lambda_{\ep, d}$ be a  subset on $X$, maximal for the property that any pair of distinct points in $\Lambda_{\ep, d}$ are distant from at least  $2\ep/\sqrt d$. Then, the union of the balls $B(x,2\ep/\sqrt d)$ centered on the points of $\Lambda_{\ep, d}$ cover $X$, and the balls $B(x,\ep/\sqrt d)$ are disjoint. Denote by $N(\Lambda_{\ep, d})$ the number of elements  $x$ of $\Lambda_{\ep, d}$ where $A(\Sigma, \mathcal L,Z(s),x, \ep/\sqrt d, c)$ happens. Then, 
by Theorem~\ref{theorem2}, 
\begin{eqnarray*}
|\Lambda_{\ep, d}|\exp\big(-\frac{c}{\ep^D}\big)& \leq & \sum_{x\in \Lambda_{\ep, d}} \mathbb P_d 
\big[
A(\Sigma, \mathcal L, Z(s),x, \ep/\sqrt d, c)
\big]\\
& = & \sum_{k=1	}^{|\Lambda_{\ep, d}|}
k\mathbb P_d ( N(\Lambda_{\ep, d})= k)\\
& \leq & \frac{1}{2}|\Lambda_{\ep, d}|e^{-\frac{c}{\ep^D}}
\mathbb P_d
\Big[
N(\Lambda_{\ep, d}) \leq \frac{1}{2}|\Lambda_{\ep, d}|e^{-\frac{c}{\ep^D}}
\Big] \\&& +
|\Lambda_{\ep, d}|\mathbb P_d
\Big[
N(\Lambda_{\ep, d}) \geq \frac{1}{2}|\Lambda_{\ep, d}|e^{-\frac{c}{\ep^D}}
\Big].
\end{eqnarray*}
Consequently, 
$\mathbb P_d
\big[
N(\Lambda_{\ep, d}) \geq \frac{1}{2}|\Lambda_{\ep, d}|e^{-\frac{c}{\ep^D}}
\big]
\geq
\frac{1}{2}
\exp\big(-\frac{c}{\ep^D}\big).$
Since 
$$\text{Vol}_{g_\omega} (X) \leq \sum_{x\in \Lambda_{\ep, d}} \text{Vol}_{g_\omega}  (x,2\ep/\sqrt d)\underset{d\to \infty}{\sim}|\Lambda_{\ep, d}|
(\frac{2\ep}{\sqrt d})^{2n} \text{Vol}_{g_0}( \mathbb B),
$$
there exists a universal $c_n>0$ 
and $d_0$  independent of $\epsilon \leq 1$ but depending on the ample probabilistic model,
such that 
$ |\Lambda_{\ep, d}| \geq  c_n \text{Vol}_{g_\omega}(X)  d^{n} \ep^{-2n}$ 
so that 
$$\mathbb P_d
\left[
N_{Lag}\big(\Sigma, \mathcal L, Z(s),x, \frac{\ep}{\sqrt d}, c\big) \geq c_n d^n  e^{-\frac{c}{\ep^D}} \text{Vol}_{g_\omega}(X)
\right]\geq \frac{1}{2}e^{-\frac{c}{\ep^D}}.$$ 
We can now absorb $c_n$ into the exponential, 
replacing $c$ by smaller positive constant. 
\end{preuve}

\begin{preuve}[ of Theorem~\ref{lagprob}]
 This is Theorem~\ref{theorem3} in the standard case and for  $r=1$.
\end{preuve}

\subsection{From  probabilistic to deterministic}

\begin{preuve}[ of Theorem~\ref{coroX}]
Theorem~\ref{coroX} is a direct consequence of Corollary~\ref{corolagrangien}   and the fact that 
the zeros of holomorphic sections  of given degree $d$ have the same diffeomorphism and symplectomorphism type, when there are equipped with the restriction of the ambient K\"ahler form $\omega$, see Proposition~\ref{nuke}.
 \end{preuve}
 \begin{remark}\label{dodo}
 As said before for projective hypersurfaces, in a parallel paper~\cite{GD}, 
 we prove the deterministic Theorem~\ref{coroX}  using the deterministic Donaldson~\cite{donaldson96} and Auroux~\cite{auroux} methods.  In the two types of proofs, we use peak sections and a lattice of mesh of order $1\sqrt d$. In both cases we prove that the Lagrangian submanifolds appear in a uniform proportion of disjoint balls centered on the vertices of the lattice.   An advantage of the Donaldson method is that it can be used for Donaldson hypersurfaces in a symplectic compact manifold $(M,\omega)$ equipped with an almost complex structure $J$. These hypersurfaces are in fact codimension 2 symplectic submanifolds which are vanishing loci  of almost holomorphic sections of high powers $L^{\otimes d}$ of a complex line bundle $L$  over $M$, where $L$ is equipped with a Hermitian metric of curvature $-2i\pi \omega$. In this general symplectic context, it is not clear which natural space of symplectic hypersurfaces can be used for probabilistic considerations. In \cite{SZA}, the authors replaced the holomorphic sections (which no longer exist in this general context) by the kernel of a certain elliptic operator acting on the bundle, which is the $\bar \partial_L$ operator if the almost complex stucture is integrable and the bundle is holomorphic. However, the vanishing locus of a section in this space is  a priori not symplectic. The deterministic proof is not easier, 
 since we also need the quantitative version of the Moser
 method given by Theorem~\ref{concrete}.
 \end{remark}

 \begin{preuve}[ of Theorem~\ref{theolag}]
  This is Theorem~\ref{coroX} in the standard case and for $r=1$.
 \end{preuve}

\subsection{Small non-contractible curves}
\begin{figure}[h]
  \centering
   \includegraphics[height=4cm]{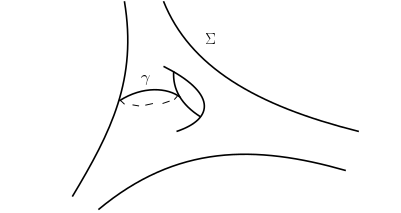}
   \caption{ A degree 3 affine complex curve $\Sigma$ in $\C^2$ with  a non-trivial loop. }\label{toto}
    \includegraphics[height=5.5cm]{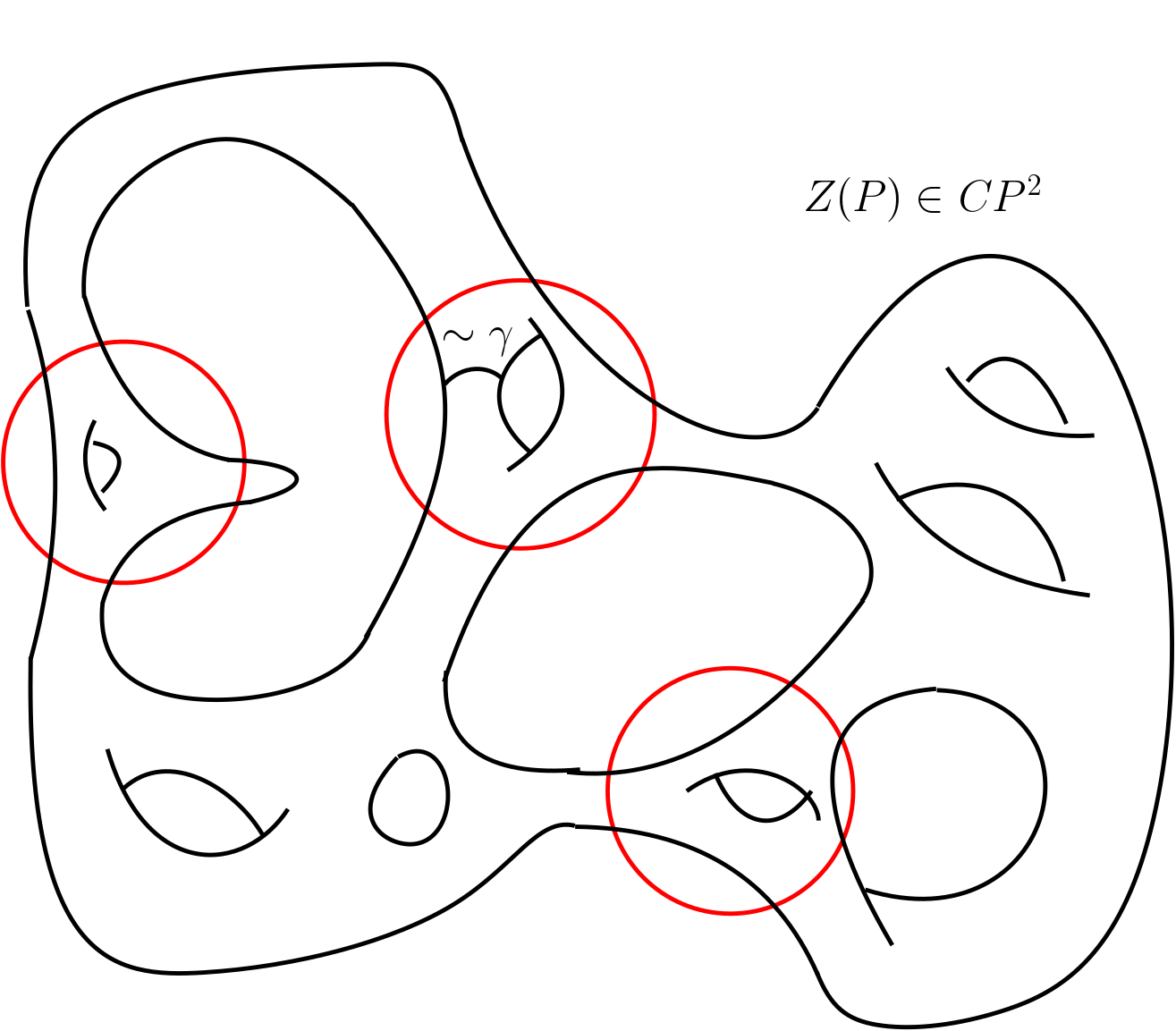}
      \caption{A (nonrealistic) degree 6 curve in $\C P^2$ and three small balls of size $1/\sqrt d$ containing the affine complex curve $\Sigma$  and the non-trivial real curve $\gamma$ of Figure~\ref{toto}.}\label{totoproj}
\end{figure}

We turn now to the proof of the Theorem~\ref{curve}
for the systoles of the random complex curves. 
\begin{preuve}[ of Theorem~\ref{curve}]
Define
$$\forall (z_1,z_2)\in \C^2, \ p( z_1,z_2) = z_1^3+z_2^3-1.$$
%$[X:Z:W]\in \C P^2$ the complex homogeneous coordinates, 
%and 
%$$\forall d, \ep>0, \ \forall (X,Z,W)\in \C^3, \ P_{\ep} (X,Z,W):= X^d p_\ep \Big( \frac{Z}{X}, \frac{W}{X} \Big).$$
%\begin{lemma}
By the genus formula applied to the 
homogenization $$P:= Z_0^3 p\big(\frac{Z_1}{Z_0},\frac{Z_2}{Z_0}, \frac{Z_3}{Z_0}\big),$$ 
$Z(P)\subset \C P^2$ is a smooth torus, so that  
for $\rho>0$ large enough, $$\tilde \Sigma:= \frac{1}\rho \big(Z(p)\cap \mathbb B(0,\rho)\big)\subset \mathbb B \subset \C^2$$ is an affine algebraic complex curve diffeomorphic to $\mathbb T^2 \setminus \cup_{i=1}^3 D_i$,
where $(D_i)_{i=1}^3$ are three disjoint disc in $\mathbb T^2$. 
Embedding $\C^2$ into $\C^n$ turns $\tilde \Sigma$ into  
an  affine algebraic complex curve $\Sigma$  in $\C^n$. 
Let $\mathcal L\subset \Sigma $ be a smooth circle which is non-trivial in 
$H_1(\Sigma,\Z )$, see Figure~\ref{toto}. Since $\mathcal L$ is a Lagrangian, by Theorem~\ref{theorem3} 
there exists at least  $d^n \text{Vol}_{g_\omega} (X)\exp(-\frac{c}{\ep^{D}})$ 
copies of $(\Sigma,\mathcal L)$ in a random curve $Z(s)$ such that any compy $\gamma_i$ of $\mathcal L$ has  intrisic diameter of the order $\ep/\sqrt d$, 
with a uniform probability given by the theorem.
The classes in $H_1(Z(s), \R)$ generated by the copies of $\gamma$ 
form an independent family. Indeed, if $\sum_i \lambda_i [\gamma_i] = 0$, 
where $(\lambda_i)_i\in \R^N$ and the $\gamma_i$ are the distinct copies of $\gamma$,
then
there exist codimension 0 surfaces with boundaries $\Sigma_1, \cdots \Sigma_{N'}$ in $Z(s)$
and $(\mu_j)_j\in \R^{N'}$ 
such that 
$\sum_i\lambda_i \gamma_i = \sum_j \mu_j \partial \Sigma_j.$ This implies
that $\partial \Sigma_j$ is a sum of distincts $\gamma_j$'s. However,
if $\gamma_i$ is one component of the boundary of 
$\Sigma_j$, then the latter   must contain the punctured torus $\tilde \Sigma$ which
contains $\gamma_i$, which implies that $\gamma_i$ bounds on the other side $\Sigma_j$,
which is a contradiction.
%if we remove a tubular neighborhood of such a real curve from $Z(s)$,
%the resulting real surface (with two boundaries) remains connected, so that if we close the two holes by two discs, 
%the Euler characteristic of the new surface is $\chi(Z(s))+1$, which implies that
%the curve was non trivial in $H_1(Z(s), \R)$ and independent from the other curves. \rk{à revoir}
%A recursion proves the affirmation. 
%\end{lemma} 
\end{preuve}
\begin{preuve}[ of Theorem~\ref{theorem1}]
Theorem~\ref{theorem1} is a particular case of Theorem~\ref{curve}, choosing $n=2$, $r=1$, $X=\C P^2$, $E=\C P^2 \times \C$, $h_E$ the Euclidean metric, $L=\mathcal O(1)$ the hyperplane bundle, $h_L$ the Fubini-Study metric and $\omega$ the Fubini-Study K\"ahler form. 
\end{preuve}
% The reason why a probabilistic method can gives  this deterministic result is the following: instead of contructing an particular hypersurface with a lot of components, like in~\cite{GDR} or~\cite{GD}, we prove in Theorem~\ref{blou} that the Lagrangian appears in the random hypersurface in a ball of size $1/\sqrt d$ centered on a fixed point, with uniform probability.

\subsection{From disjoint to homologically non-trivial}

\begin{preuve}[ of Corollary~\ref{corolagX}]
The first assertion is a direct consequence of the classical Lemma~\ref{totally} below, remembering that Lagrangian submanifold are totally real 
 for any almost complex stucture tamed by the symplectic form $\omega$. The second assertion was explained in Remark~\ref{remarkseveral}.
 \end{preuve}
 \begin{lemma}\label{totally}
 Let $\mathcal L\subset (X,J)$ be any closed oriented smooth totally real dimension $n$ submanifold in an almost complex manifold $X$ of dimension $2n$. Then,
 $$ [\mathcal L]\cdot [\mathcal L]= \chi(\mathcal L),$$
 where $[\mathcal L]\in H_n(X,\Z)$ and $\chi(\mathcal L)$ denotes the Euler characteristic of $\mathcal L$. 
 If $\mathcal L_1, \cdots, \mathcal L_N$ is a family of disjoint totally real submanifolds of $X$ with nonvanishing
 Euler characteristic, then the family made of their classes $[\mathcal L_1], \cdots, [\mathcal L_N]$ 
 in $H_{n}(X,\R)$ is independent. 
 \end{lemma}
 \begin{preuve}
 For a closed totally real $\mathcal L\subset X$, if $h$ is any metric, then $JT\mathcal L \sim NL$,
 where $N\mathcal L$ is the normal bundle over $\mathcal L$. Then
 $\chi(\mathcal L) = \int_{\mathcal L} e(T\mathcal L) = \int_L e(N\mathcal L)
 $ which equals $[\mathcal L|\cdot [\mathcal L] $.
 For the second assertion, 
if $\sum_{i=k}^m a_i [\mathcal L_k]= 0$ in $H_n(X, \R)$, where $\mathcal L_1, \cdots , \mathcal L_m$ are pairwise disjoint totally real submanifolds, then for every $j$, intersecting by $[\mathcal L_j]$ gives $a_j [\mathcal L_j]\cdot [\mathcal L_j]=0$ so that in our case, $a_j= 0$.
\end{preuve}

\subsection{From smooth to algebraic}
For the proof of Corollary~\ref{corolagrangien} we will use the
classical theorem of H. Seifert:
\begin{theorem}{\cite{seifert}}\label{seifert}
Let $n\geq 2$, $1\leq r \leq n$ and $\mathcal L\subset \R^n$ be any compact smooth $(n-r)$-submanifold with trivial normal bundle. Then, there exists a real
polynomial map $p:= (p_1, \cdots, p_r): \R^n \to \R^r$ and 
a diffeotopy of $\R^n$ sending $\mathcal L$ onto some connected components of  $Z_\R(p).$ The diffeotopy can be chosen as $C^1$-close to the identity map as we want. 
\end{theorem}
It is not known which hypersurfaces are diffeotopic to algebraic ones, see~\cite[Remark 14.1.1]{bochnak2013real}.

\begin{preuve}[ of Corollary~\ref{corolagrangien}]
By Theorem~\ref{seifert}, there exists a regular real polynomial  map 
$p = (p_1, \cdots, p_r): \R^n \to \R^p$ of maximal degree $d(p):= \max_i \deg p_i$  such that $Z_{\R}(p):= Z(p)\cap \R^n$ has a compact component $\mathcal L'$ or a set of components diffeomorphic to $\mathcal L$. If $\mathcal L$ is algebraic, we can choose $p$ such that $Z_{\R}(p)= \mathcal L$. By a comprehensible abuse of notation, we keep the notation $\mathcal L$ for $\mathcal L$. After perturbation, we can assume that $p$, when considered as  defined on $\C^n$, is regular, too. Then, $Z_\R(p)$ is  a Lagrangian submanifold
of its complex vanishing locus $\Sigma := Z(p)$ equipped by the the restriction of the standard K\"ahler form $\omega_0$. 
For a large enough $\rho>0$, 
$\rho\mathbb B $ contains $\mathcal L$.
We rescale the polynomial by $1/\rho$ and keep the notation $p$, so that $Z_\R (p) \cap \mathbb B $ contains $\mathcal L$. 
Then Corollary~\ref{corolagrangien} is a consequence of Theorem~\ref{theorem3} applied to the couple $(\Sigma\cap \mathbb B, \mathcal L)$.
\end{preuve}
% % % % % % % % % % % % % %
\section{Quantitative deformations}\label{quant}
In this section we introduce and prove
deterministic lemmas and propositions
which quantify how much a given 
specific geometrical situation 
can be perturbed keeping its specificity. The first part concerns the topology, the second part
being Lagrangian. 
% % % % %
\subsection{Preserving the topology}
The next proposition
is a quantitative and deterministic version of the barrier method for  functions. 
We need first a notation. 
For any linear mapping $A\in \mathcal L (\R^m, \R^p)$, where $1\leq p\leq m$ are integers,   define
\begin{equation}\label{T}
 T(A) := \inf_{|w|=1} |A^* w|,
 \end{equation}
where $| \cdot |$ denotes the standard Euclidean norm. 
In the sequel we will use the following simple properties: for any $A\in \mathcal L(\R^m, \R^p), $
\begin{itemize}
\item $T(A)>0$ if and only if $A$ is onto;
\item $T(A)\leq \|A\|$, where $
\|A\|:= \sup_{|v|=1} |Av|$;
\item $\|(AA^*)^{-1}\|\leq T(A)^{-2}$;
\item if $p=1$, then $T(A) = \|A\|$;
\item for any $ B\in \mathcal L(\R^m, \R^p), \ T(A+B)\geq T(A) - \|B\|$. 
\end{itemize}
The following proposition provides quantitative
estimates for the perturbation  of 
a vanishing locus on $2\mathbb B$. It differs from 
\cite[Proposition 3.4]{GWL} in two ways. 
First, it allows the vanishing locus to cross the boundary 
of the ball. Second, it explicits quantitatively 
the existence of a diffeomorphism sending
the vanishing locus to its perturbation. 
We need indeed to understand how a Lagrangian submanifold
of the locus can be moved into another Lagrangian submanifold
of the perturbed locus. For  this,  we give quantitative estimates of the difference between 
the  diffeomorphism and the identity. 
For $\eta>0$ and $f,g: 2\mathbb B \to \R^p$ $C^k$ mappings,
define
\begin{equation}\label{cocorico}
\forall 0\leq j\leq k-1, \ c_j(\eta, f,g) := \frac{1}{\eta^{2(j+1)} }\|f\|^{2j+1}_{C^{j+1}(2\mathbb B)}\|g\|_{C^{j}(2\mathbb B)}.
\end{equation}
Note that $c_j$ is an homogeneous function of degree 0; this will be crucial for probabilistic estimates, see~(\ref{chouchou}) below. 
\begin{proposition}\label{morse} Let $m\geq 1$, $1\leq p\leq m$ and $k\geq 3$ be integers,   $\eta>0,$
and $ f,g : 2\overline{\mathbb B}\subset \R^m  \to \R^p$ be two  $C^k$ maps, such that   $ \|g\|_{C^1(2\mathbb B)} \leq \eta/8$,  
$ c_0(\eta, f, g)\leq  1/8$ and
$$
 \forall x\in 2\mathbb B, \ |f(x)|<\eta\Rightarrow T\big(df(x)\big)>\eta.$$
\begin{enumerate}
\item \label{morse1}
Then, there exists a $1-$parameter family $(\phi_t)_{t\in [0,1]}$ of 
diffeomorphisms with support in  $2\mathbb B$ such that  
 $$\forall t\in [0,1], \ \big(Z(f), \overline{\mathbb B }\big)
 \sim_{\phi_t}
 \big(Z(f+tg) ,\phi_t\big(\overline{\mathbb B}\big)\big)$$ 
 with $\ Z(f+tg)\cap \frac{1}{2}\mathbb B\subset  \phi_t(Z(f)\cap \mathbb B) \subset Z(f+tg)\cap \frac{3}{2}\mathbb B$,
  $(x,t)\mapsto \phi_t(x)$ is $C^{k-1}$ and 
 \begin{eqnarray}\label{phi}
 \forall t\in [0,1],  \ 
 \|\phi_t- \id\|_{C^0(2\mathbb B)}&\leq &tc_0 (\eta, f, g).
   \end{eqnarray} 
 \item   
 \label{morse2}
 Let $j=1,2$ and  $C>1$ such that $c_j(\eta, f,g)\leq C$. Then there exists $C'$ 
 depending only on $C$
 such that 
    \begin{eqnarray}\label{phi1}
  \forall t\in [0,1],  \    \|\phi_t- \id\|_{C^j(2\mathbb B)}&\leq & C'tc_j(\eta, f, g).
      \end{eqnarray} 
\end{enumerate}
 \end{proposition}
 % % % % % % % % % % % % % % % %
 In the proof of the main probabilistic Theorem~\ref{theorem2} below, the different  estimates for the various norms of $\phi_1-\id$ in Proposition~\ref{morse}
will be used in different ways:
 \begin{itemize}
\item   a small $C^0$ norm will imply
that a Lagrangian submanifold of $Z(f)$ in $\frac{1}{2}\mathbb B$
will be send by $\phi_1$ in $Z(f+g)$ into a submanifold
of $\mathbb B$;
\item  a small $C^1$ norm implies that $\phi_1$ is close to be symplectic, so that the image of the Lagrangian is close to be Lagrangian and can be perturbed into a genuine Lagrangian submanifold of $Z(f+g)$, see Theorem~\ref{concrete}; 
\item the bound for the $C^2$ norm will be used to estimate the
intrisic metric on the perturbation of $Z(f)$ on $Z(f+g)$, in order ton obtain the estimates for diameters. 
\end{itemize}
% % % % % % % % % % % % % % % % % % % % % % % % % %
 \begin{preuve}[ of Proposition~\ref{morse}]
 Define for any $t\in [0,1]$, 
$  f_t :=  f+tg$. 
We first prove that 
\begin{equation}\label{transy}
\forall (x,t)\in 2\mathbb B\times [0,1], \ |f_t(x)|< \eta/2\Rightarrow T\big(df_t(x)\big)>\eta/2.
\end{equation}
Indeed, 
$|f_t(x)|< \eta/2$ implies $|f(x)|<\eta$ since $|g(x)|<\eta/2$, so that $T(df(x) )>\eta$ by hypothesis, and since 
 $df_t =  df + tdg $ and $\|dg(x)\|<\eta/2$, then
 $T(df_t(x))>\eta/2$. In particular, for all $t\in [0,1]$,
$Z(f_t)$ is a $C^{k-1}$ codimension $p$ submanifold of $2\mathbb B$.
For any $t\in [0,1]$, $\beta> 0$, let 
 $V_t(\beta):= \{x\in 2\mathbb B, |f_t(x)|\leq \beta\}.$
Then, by hypothesis on $g$, 
\begin{equation}\label{inclusions}
\forall t\in [0,1], \ Z(f_t)\subset V_0\big( \eta/8\big)\subset V_0\big(\eta/4\big)\subset V_t\big(\eta/2\big).
\end{equation}
For all $(x,t)\in V_0\big(\eta/4\big)\times [0,1]$ define 
$  X_t(x)\in \R^m$ to be the projection of the origin 
onto the $(m-p)$-plane $$df_t(x)^{(-1)}(\{-\partial_t f_t(x)\})\subset \R^m,$$ 
which is well defined by  (\ref{inclusions}) and (\ref{transy}).
Note that  $X(x,t) = \Phi (df_t(x), g(x))$ 
  where $\Phi$ is defined in Lemma~\ref{proj} of the annex. This Lemma~\ref{proj} shows that $\Phi$ is a smooth mapping where the first variable is onto, 
  so that $X$ is $C^{k-1}$ where it is defined. 
 Let $\chi : \R\to [0,1]$ be a smooth cut-off function satisfying $\chi_{|(-\infty ,1/4]}= 1$ and $\chi_{|[1/2,1]}= 0$, 
 and define on $2\mathbb B$ the family of vector fields $$
 \forall (x,t)\in 2\mathbb B\times [0,1], \ 
 \tilde X_t(x) := \chi \Big(\frac{2}{\eta}|f(x)|\Big) \chi \Big(\frac{|x|-1}2 \Big) X_t(x).$$
 Then, $\tilde X_t $ is $C^{k-1}$ in $(t,x)$, for any $t\in [0,1]$ $ \tilde X_t= X_t$ over $V_0\big(\eta/8\big)\cap \frac{3}{2}\mathbb B$, and $\tilde X_t = 0$ 
 on $\big(V_0(\eta/4)\big)^c$ and on $\partial (2\mathbb B)$.  
 %and $Y$ is unchanged on a neighourhood of $Z(F)$.  
 Now define
 $(\phi_t)_{t\in [0,1]}$ the family of diffeomorphisms
 generated by $(\tilde X_t)_{t\in [0,1]}$ on $2\mathbb B $,
 that is 
$$\forall (x,t)\in 2\mathbb B\times [0,1], \ \partial_t \phi_t (x) = \tilde X_t\big(\phi_t(x)\big), \ \phi_0 = \id.$$
Note that $(x,t)\mapsto \phi_t(x)$ is $C^{k-1}$.
By construction,  $\phi_t $ can be extend smoothly as the identity outside $2\mathbb B$. Since the $C^0$ norm of $\tilde X$ is bounded by the one of $X$, 
  by Lemma~\ref{majphi} and Lemma~\ref{proj}, 
  $$ \forall t\in [0,1], \ \|\phi_t-\id\|_{C^0(2\mathbb B)}\leq \frac{4t}{\eta^2}\|df\|_{C^0(2\mathbb B)}\|g\|_{C^0(2\mathbb B)}\leq  1/2$$
by hypothesis, so that 
$\forall t\in [0,1], \  \frac{1}{2}\mathbb B\subset \phi_t(\mathbb B) \subset \frac{3}{2}\mathbb B.$
Moreover, for any $(x,t)$ such that $\phi_t(x) \in V_0\big(\eta/8\big)\cap \frac{3}{2}\mathbb B$, $$ \partial_t \big(f_t( \phi_t(x))\big) = 
g( \phi_t(x)) + df_t( \phi_t (x)(X(\phi_t(x), t) = 0.$$
By an open-closed argument and the inclusions~(\ref{inclusions}), this condition is satisfied if $x\in \mathbb B$.
Consequently,
$$\forall t\in [0,1], \ Z(f_t)\cap \frac{1}{2}\mathbb B\subset  \phi_t(Z(f)\cap \mathbb B) \subset Z(f_t)\cap \frac{3}{2}\mathbb B$$
and the first assertion \ref{morse1}. of the proposition is proved. 

Now, since 
$
  \forall t\in [0,1], \  dX_t= d\Phi \big(df_t(x), g(x)\big) (d^2 f_t , dg(x)),
 $
 Lemma~\ref{proj} gives 
 \begin{eqnarray*}
  \max_{t\in [0,1]} \|dX_t\|_{C^0(2\mathbb B)} &\leq  &
 \frac{16}{\eta^4}\|g\|_{C^0(2\mathbb B)} 
 \big(2\|df\|_{C^0}^2 +\eta^2/4)\|d^2f\|_{C^0}
 \big)+ \frac{4}{\eta^2}\|df\|_{C^0} \|dg\|_{C^0}\\
 &\leq & 
Kc_1(\eta, f, g),
 \end{eqnarray*}
 where $K$ is a universal constant.
Moreover,
$ d\tilde X_t = \frac{\nabla |f|}{\eta} \chi' X+\chi dX,$
so that $$\|d\tilde X_t\|_{C^0(2\mathbb B)}\leq K'\big(\frac{\|df\|_{C^0}}{\eta}c_0(\eta, f, g)+ c_1(\eta, f, g)\big)
\leq K'' c_1(\eta, f, g),$$
 where $K', K'' $ depend only on $\chi$. Consequently, by Lemma~\ref{majphi}, there exists
$C'$ depending only on  $C$, such that 
$$ \|\phi_t - \id\|_{C^1(2\mathbb B)} \leq C'tc_1(\eta, f, g).$$
This proves assertion~\ref{morse2}. ~for $j=1$.

For $j=2$ in assertion~\ref{morse2}., we compute
   \begin{eqnarray*}
  d^2_xX_t(x,t)&=& d^2\Phi \big(df_t(x), g(x)\big) (d^2 f_t , dg)^2
  + d\Phi \big(df_t(x), g(x)\big) (d^3 f_t , d^2g(x)).
  \end{eqnarray*}
  so that
  by Lemma~\ref{proj},
 \begin{eqnarray*}
   \|d^2_xX_t\|&\leq & 14\|df\|^3\eta^{-6} |g| \|d^2f\|^2 + 6\eta^{-4} \|df\|^2 \|d^2f \| \|dg\|+\\
   && 3\eta^{-4} \|df\|^2 |g| \|d^3 f\|+ \eta^{-2}\|df\| \|d^2g\|\\
   & \leq & 24c_2(\eta, f,g).
   \end{eqnarray*} 
A similar estimate for $d^2\tilde X$ and Lemma~\ref{majphi} imply
$$ \forall t\in [0,1], \ \|\phi_t-\id\|_{C^2(2\mathbb B)}\leq 
tC''c_2(\eta, f, g),$$
where $C''$ is a constant depending only on $C$. 
 \end{preuve}

\subsection{Preserving the Lagrangianity}
The main goal of this paragraph is to prove the 
technical Proposition~\ref{defo} below. The latter
asserts that, in a quantitative way,  
\begin{itemize}
\item 
if some compact Lagrangian submanifold $\mathcal L$ 
lies inside a compact  symplectic submanifold $\Sigma$ of a symplectic manifold $(M,\omega)$, 
\item if $\Sigma$ is perturbed into $\phi(\Sigma)$ by a diffeomorphism $\phi$ close to the identity, 
\item and if $\omega$ is exact and perturbed by a small 2-form $d\mu$, 
\end{itemize}
then there exists a perturbation $\mathcal L'\subset \phi(\Sigma)$ of $\mathcal L$ which is Lagrangian for the restriction of the perturbed form. 
Since we think that
this quantitative proposition has its own interest, we provide a general statement and a proof for symplectic manifolds. However, in this paper we will apply it in the simple case where the ambient manifold is the unit ball of the standard symplectic space $(\R^{2n}, \omega_0)$, see Theorem~\ref{concrete} below. 
Before the statement of Proposition~\ref{defo}, we need some definitions:
\begin{definition}
Let $(M, h)$  be  a smooth Riemannian manifold, possibly with boundary.
\begin{itemize}
\item For any  pair of continuous maps $f, g : M\to M$
define
$ d(f,g) := \sup_{x\in M} d(f(x), g(x))$, where $d$ is the distance associated to $h$. 
\item  For any $k\geq 0$ and any  $C^k$ vector field $X$ on $M$, define
$ N_k(X,M)= \sup_{x\in M, 0\leq p\leq k} \|\nabla^p X\|$
and similarly $N_k(\alpha, M) $ for any $C^k$ form $\alpha$ on $M$.
Here, $\nabla$ denotes the Levi-Civita connection associated to $h$. 
\item  For any submanifold $\mathcal L\subset M$, define $\diam_M(\mathcal L):= \max_{p,q\in \mathcal L} d(p,q).$
\item For any 2-form $\omega$ defined on a neighborhood of an open $U$ set of a manifold $M$ equipped with a metric $h$, 
let 
$S(\omega,U) : = \inf_{x\in U, X\in T_x M, |X|=1 } \sup_{Y\in T_xM, |Y|=1} |\omega(X,Y)|.$
\end{itemize}
\end{definition} 
Note that: 
\begin{itemize}
\item $\diam_{\mathcal L}(\mathcal L)$ is the intrisic diameter of $\mathcal L$. Note also that if $\mathcal L$ is a circle, then its length is bounded by its intrisic diameter.
\item 
If $U$ is relatively compact, then 
$\omega$ is symplectic over $U$ if and only if $S(\omega, U)>0$; 
\item if $\omega_0$ denotes the standard symplectic form on $\R^{2n}$, that is $\omega_0= \sum_{i=1}^n dx_i\wedge dy_i$, then $S(\omega_0, \R^{2n}) = 1$;
\item 
for any pairs of 2-forms $\omega$ and $\omega'$,
\begin{equation}\label{Caligula}
S(\omega+\omega', U) \geq S(\omega, U) - \|\omega'\|;
\end{equation}
\item if $X$ is a vector field on $U$, $i_X \omega = \lambda$ and $\omega$ is symplectic, then 
\begin{equation}\label{Neron}
\|X\|_{C^0(U)} \leq \frac{1}{S(\omega, U)}\|\lambda\|_{C^0(U)};
\end{equation}
\item if $f: U\subset \C^n \to \C$ is holomorphic and vanishes transversally, then $S(\omega_{0|TZ(f)}, U) = 1$. 
\end{itemize}
\begin{proposition}\label{defo}
Let  $1\leq r\leq n$ be integers, $(M,\omega,h)$ a smooth symplectic $2n-$manifold
equipped with a metric, 
$U\subset V\subset W\subset Y$ four relatively compact open sets such that $\overline{U}\subset V$, $\overline V\subset W$, $\overline W\subset Y$, and assume that there exists $\lambda $ a smooth 1-form on $\overline{Y}$ such that $\omega_{|\overline{Y}} = d\lambda$. Let 
$\Sigma \subset \overline{Y}$ be a compact smooth codimension $2r$ submanifold,   symplectic for $\omega_{|T\Sigma}$,  
$\mathcal L$ be a compact smooth Lagrangian submanifold of $(\Sigma\cap U, \omega_{|T\Sigma})$, 
$\phi : {Y} \to {Y} $ be a smooth diffeomorphism with support in $Y$, and $\mu$  be a smooth 1- form on $\overline{W}$ satisfying
 \begin{eqnarray*}
d(\phi, \id) &\leq &\text{dist}(V,\partial W),\\
 \|\phi^*(\lambda +\mu) -\lambda\|_{C^0(W)} 
 &\leq & \frac{1}{2}S(\omega_{|T\Sigma}, W\cap \Sigma) \text{dist}(U,\partial V) \\
\text{ and } \|\phi^*(d\lambda +d\mu) -d\lambda\|_{C^0(W)} &\leq & \frac{1}{2}S(\omega_{|T\Sigma}, W\cap \Sigma).
 \end{eqnarray*} 
 \begin{enumerate}
 \item \label{defo1} Then, there exists $\mathcal L' $ 
 a compact smooth Lagrangian submanifold of $\big( \phi (\Sigma)\cap W, (\omega+d\mu)_{|\phi(\Sigma)}\big), $ such that
 $\big(\mathcal L, \Sigma\cap V)\sim_{\phi}\big(\mathcal L', \phi(\Sigma\cap V)\big).$
 % % % % %
 \item \label{defo2}If furthermore  $d(\phi, \id) \leq    \frac{1}{8}\diam_M(\mathcal L)$ and
  \begin{eqnarray*}
   \|\phi^*(\lambda +\mu) -\lambda\|_{C^0(W)} 
     \leq  \frac{1}{16}S(\omega_{|T\Sigma}, W\cap \Sigma)  \diam_M(\mathcal L), 
  \end{eqnarray*} 
 then 
$ \frac{1}{2}\diam_M (\mathcal L)\leq \diam_{\mathcal L'} (\mathcal L').$
\item \label{defo3} Let $C>1$.  If furthermore
$$\max \Big(S(\omega_{|T\Sigma}, W\cap \Sigma)^{-1}, N_1\big(\phi^*(\lambda +\mu)  -\lambda,W\big), N_1(\omega, W), \|d\phi\|_{C^0(W)} \Big) \leq  C,$$
 then there exists $C'>0$ depending only on $C$, on the pair $(V,W)$ and on the $C^1$ norm of $h$ over $W$ such that 
$ \diam_{\mathcal L'} (\mathcal L') \leq C'\diam_{\mathcal L} (\mathcal L).$
 \end{enumerate}
\end{proposition}
In the proof of Theorem~\ref{theorem2}, where we prove that a given affine complex hypersurface $\Sigma$ with a Lagrangian submanifold  $\mathcal L$ appears with uniform probability in a sequence of small balls, we will need  Proposition~\ref{defo}
applied to the concrete context of Proposition~\ref{morse}, where $\Sigma$ is the vanishing
locus of a holomorphic function $f$, $\phi$  is a diffeomorphism
sending $Z(f)$ onto the perturbed submanifold $Z(f+g)$ and $\omega$ the K\"ahler form viewed in the chart on the standard ball.
The next Theorem~\ref{concrete} below synthesizes these two propositions for this goal: it asserts, in a quantiative way,  that if $\mathcal L$ is a compact Lagrangian of a vanishing locus $Z(f)$ which is symplectic
for the restriction of the standard form inside the standard ball, as it is the case for the real part of a complex hypersurface defined by a real polynomial, 
and if $g$ is a small perturbing function, then there exists a perturbation $\mathcal L'$ of $\mathcal L$ which is a Lagrangian submanifold of 
$Z(f+g)$ equipped with the restriction of a  perturbation $\omega_0+d\mu$ of the standard form. 
% % % % % % % % % % % % % % % % % % % % % % % %
\begin{theorem}\label{concrete}
 Let $n\geq 1$ and $1\leq r\leq n$ be integers,   $\eta>0$,  
and $ f,g : 2\overline{\mathbb B}\subset \R^{2n}  \to \R^{2r}$ be two  smooth maps, such that $\|g\|_{C^1(2\mathbb B)} \leq \eta/8,$ and
$$
 \forall x\in 2\mathbb B, \ |f(x)|<\eta \Rightarrow T\big(df(x)\big)>\eta.$$
Let $\omega$ be a smooth symplectic form on $2\overline{\mathbb B}$ and  $\mu$ be a smooth 1-form on $2\overline{\mathbb B}$
satisfying  $\omega= \omega_0+d\mu$, with
\begin{eqnarray}\label{kikou}
\max\big(c_0(\eta, f,g), c_1(\eta, f, g), 
\|\mu\|_{C^0(2\mathbb B)},\|d\mu\|_{C^0(2\mathbb B)}\big) &\leq  &
\frac{1}{16}.
\end{eqnarray}
Let $\mathcal L$ be a compact smooth Lagrangian submanifold of $\big( Z(f)\cap \frac{1}{2}\mathbb B, \omega_{0|Z(f)}\big)$.  
Then, 
\begin{enumerate}
\item\label{con1} there exists a smooth ball $B$ satisfying $\frac{1}{2}\mathbb B \subset B\subset \frac{3}{2}\mathbb B$ and $ \mathcal L'$ a compact smooth Lagrangian submanifold of $\big(Z(f+g)\cap B, \omega_{|Z(f+g)}\big)$,
satisfying 
$ (\mathcal L, Z(f)\cap \mathbb B)  \sim_{diff} (\mathcal L', Z(f+g)\cap B)$.
% % % % % %
\item\label{con2} If furthermore 
\begin{eqnarray}\label{cloclo}
\max\big(c_0(\eta, f,g), c_1(\eta, f, g), 
\|\mu\|_{C^0(2\mathbb B)}\big) &\leq  &
\frac{1}{16}  \diam_{Z(f)}(\mathcal L),
\end{eqnarray}
then
$ \frac{1}{2}\diam_{Z(f)} (\mathcal L)\leq \diam_{\mathcal L'} (\mathcal L').$
% % % % % %
 \item \label{con3} Let $C>1$ be such that, furthermore, 
\begin{equation}\label{aragon}\max\big(  c_2(\eta, f,g), N_1(\mu, 2\mathbb B)\big)\leq C.
 \end{equation}
 Then there exists $C''>1$ depending only on $C$,  such that
$ \diam_{\mathcal L'} (\mathcal L') \leq C'' \diam_{\mathcal L} (\mathcal L).$
\end{enumerate}
\end{theorem}
% % % % % % % % % % % % % % % % % % % % % % % % % % %
The various estimates for the diameters concern the restriction of the standard metric $g_0$ on $\R^{2n}$. We postpone the proof of Proposition~\ref{defo} and prove the theorem now, which is the consequence
of the latter Proposition~\ref{defo} and the former Proposition~\ref{morse}.
\begin{preuve}[ of Theorem~\ref{concrete}]
By the two first assertions~\ref{morse1}. and \ref{morse2}. of Proposition~\ref{morse}, there exists a family of diffeomorphisms 
$(\phi_t)_{t\in [0,1]} :2\mathbb B \to  2\mathbb B$ with compact support and a universal constant $K'\geq 1$ such that,
writing $\phi=\phi_1$, 
$$ d(\phi, \id)= \|\phi-\id\|_{C^0} \leq c_0(\eta, f,g)\leq \frac{1}{2} \text{ and } \|d\phi- \id\|_{C^0}\leq K'c_1(\eta, f,g), $$
and $(Z(f), \mathbb B) \sim_{\phi} (Z(f+g), \phi(\mathbb B)$ with $$Z(f+g)\cap \mathbb B\subset \phi(Z(f)\cap \mathbb B)
\subset Z(f+g)\cap \frac{3}{2}\mathbb B.$$
%Since $\|d\phi\|\leq 1+K'/16$, there exists a universal constant $K''>0$ such that $$\diam_{Z(f+g)} Z(f+g)\cap \phi(\mathbb B)\leq K''\diam_{Z(f)} Z(f)\cap \mathbb B.$$ 
 
 Let $ \lambda_0 := \sum_{i=1}^n 
 x_idy_i$ be the standard Liouville form, which satisfies $d\lambda_0 = \omega_0$. Note that  for any $x\in \R^{2n}$, $\|\lambda_0(x)\|\leq |x|.$ 
Then, using that $S(\omega_{0|TZ(f)},2\mathbb B)=1$, 
\begin{eqnarray}
    \|\phi^*(\lambda_0 +\mu) -\lambda_0\|_{C^0(2\mathbb B)}
&\leq &
\|\phi-\id\|_{C^0}+
\|\lambda_0\|_{C^0}\|d\phi-\id\|_{C^0}+\|d\phi\|_{C^0} \|\mu\|_{C^0}\nonumber \\
%&\leq & \frac{1}{\eta^2 }\|f\|_{C^1}\|g\|_{C^0} +\frac{2K}{\eta^4}\|f\|^3_{C^2}\|g\|_{C^1}+\big(1+\frac{2K}{\eta^4}\|f\|^3_{C^2}\|g\|_{C^1}\big) \|\mu\|_{C^0}\nonumber\\
&\leq &  
c_0(\eta, f, g) +2c_1(\eta, f,g) +(1+c_1) \|\mu\|_{C^0} \nonumber\\
&\leq &  \label{forme} 5\max(c_0, c_1, \|\mu\|_{C^0})\\
&\leq & \frac{1}{2}S\big(\omega_{0|TZ(f)}, 2\mathbb B\cap Z(f)\big)\text{ dist}(\mathbb B, 2\mathbb B)  \nonumber
%&\leq & 3tC' \|g\|_{C^1} + (1+tC')\|\mu\|_{C^0}\\
%&\leq & 3C'( \|g\|_{C^1(2\mathbb B)} +\|\mu\|_{C^0(2\mathbb B})\\
%&\leq & \frac{1}{16C^2} \min \big(\text{dist}(\frac{1}{2}\mathbb B,\partial \mathbb B), \diam_{\R^{2n}}(\mathcal L)\big)
\end{eqnarray}
by (\ref{kikou}).
Similarly, 
\begin{eqnarray*}
  \|\phi^*(d\lambda_0 +d\mu) -d\lambda_0\|_{C^0(2\mathbb B)}
&\leq &
\|d\phi-\id\|_{C^0}^2+
2\|d\phi-\id\|_{C^0}+\|d\phi\|^2_{C^0} \|d\mu\|_{C^0}\\
&\leq & c_1^2+2c_1+(1+c_1)^2 \|d\mu\|_{C^0}\\
&\leq & 
\frac{1}{2}S\big(\omega_{0|TZ(f)}, 2\mathbb B\cap Z(f)\big)
\end{eqnarray*}
again by (\ref{kikou}).
 By assertion \ref{defo1}. of Proposition~\ref{defo} applied to  $Y=2\mathbb B$, $W= \frac{3}{2}\mathbb B$, $V= \mathbb B$, $U= \frac{1}{2}\mathbb B$, $\Sigma= Z(f)$, 
and  $\omega = \omega_0$,  there exists a Lagrangian submanifold 
$\mathcal L'$ 
of $\big( Z(f+g)\cap \phi(\mathbb B), (\omega_0+d\mu)_{|TZ(f+g)}\big), $ such that
$$\big(\mathcal L, Z(f)\cap \mathbb B)\sim_{\phi} \big(\mathcal L', Z(f+g)\cap \phi(\mathbb B) \big).$$
If $B:= \phi(\mathbb B)$, then  $\frac{1}{2}\mathbb B\subset B\subset \frac{3}{2}\mathbb B$. Hence, the first assertion of the theorem is proved. 
% % % % % % %
%
If furthermore~(\ref{cloclo}) is satisfied,  using~(\ref{forme}), the hypotheses of assertion \ref{defo2}. of Proposition~\ref{defo} are satisfied, so that 
$\frac{1}{2}\diam_{Z(f)} (\mathcal L)\leq \diam (\mathcal L') $.
This proves the second assertion. 
% % % % % %
%
Now, if (\ref{aragon}) is satisfied, by assertion \ref{morse2}. of Proposition~\ref{morse},  there exists $C''$ depending only on $C$ such that
$\|d^2\phi\|_{C^0(2\mathbb B)}\leq C''$. 
This implies  that  
there is a universal constant $K'''$ and a constant $C'''$ depending only on $C$  such that 
$$N_1 \big[\phi^*(d\lambda_0 +d\mu) -d\lambda_0), 2\mathbb B\big]\leq K'''
\big(\|d\phi\|_{C^0}N_1(d\mu, 2\mathbb B)+\|d^2\phi\|_{C^0}\|d\mu\|_{C^0}\big)\leq C''.$$
Consequently,
$$\max \Big(S\big(\omega_{0|TZ(f)}, \frac{3}{2}\mathbb B\cap Z(f)\big)^{-1},  N_1 \big(\phi^*(d\lambda_0 +d\mu) -d\lambda_0, \frac{3}{2}\mathbb B\big), N_1(\omega_0, 2\mathbb B)\Big) \leq C''',$$ 
where $C'''$ depends only on $C$. 
We can now apply assertion~\ref{defo3}. of Proposition~\ref{defo}: 
there exists a constant $C''''>0$ depending only on $C$ such that
$  \diam_{\mathcal L'} (\mathcal L') \leq C'\diam_{\mathcal L} (\mathcal L).$
\end{preuve}

The main steps for proving Proposition~\ref{defo} are the following:
\begin{itemize}
\item Recall that in Proposition~\ref{defo}, a symplectic submanifold $\Sigma$ is deformed by a diffeomorphism $\phi$ into $\Sigma'$, and the ambient symplectic form $\omega$ is deformed into $\omega+d\mu$. 
\item The restriction of the perturbed form on $\Sigma'$ can be viewed as $\Omega'= \phi^*(\omega+d\mu)_{|T\Sigma}$ on $\Sigma$.  Proposition~\ref{moser} 
constructs, in a general setting, an isotopy of local diffeomorphisms $(\psi_t)_t$ on $\Sigma$, such that $\psi_1$ is a symplectomorphism between $\Omega'$ and a given symplectic form $\Omega$, which is $\Omega= \omega_{|T\Sigma}$ in our case,  with an explicit control of $\psi_1-\id$ depending on $\Omega-\Omega'$ and its primitive. 
\item Corollary~\ref{perturbation} applies this intrisic Proposition~\ref{moser} to the relative situation of Proposition~\ref{defo}, and transfers the latter control to controls depending on $\phi-\id$ and the perturbation of the ambient symplectic form. 
\item The proof of Proposition~\ref{defo} consists in applying this corollary to the deformation of the Lagrangian submanifold. 
\end{itemize} 
{\bf Moser trick.} 
The next Proposition is a quantitative  version of the Moser's trick.
\begin{proposition}\label{moser} Let $(\Sigma,\Omega, H)$ be a smooth symplectic manifold, possibly with boundary, equipped with a metric $H$, 
and $\mathcal U, \mathcal V,  \mathcal W$ be three relatively compact open sets in $\Sigma$ such that 
$\overline {\mathcal U}\subset \mathcal V$ and 
$\overline{ \mathcal V}\subset \mathcal W$.
Let $\nu$ be a smooth 1-form on $\overline{\mathcal W}$
satisfying  $$\|\nu\|_{C^0(\mathcal W)}\leq \frac{1}{2}S(\Omega, \mathcal W))\text{dist} (\mathcal U, \partial \mathcal V) \text{ and }\|d\nu\|_{C^0(\mathcal W)}\leq \frac{1}{2}S(\Omega, \mathcal W).$$ 
\begin{enumerate}
\item \label{moser1} Then, there exists 
a smooth family of diffeomorphisms $(\psi_t)_{t\in [0,1]} : \mathcal W\to \mathcal W$ with support in $\mathcal W$ such that
$$\forall t\in [0,1], \ \psi_t^*\big(\Omega+  td\nu\big) = \Omega \text{ on }\mathcal U, \  \psi_t(\mathcal U)\subset \mathcal V, \text{ and }  d(\psi_t, \id)\leq \frac{2t}{S(\Omega, \mathcal W)}\|\nu\|_{C^0(\mathcal W)}.$$
% % % % % % % % %
\item \label{moser2} Let $ C>1$ and 
assume that 
\begin{equation}\label{stirner}
\max \big( S(\Omega, \mathcal W)^{-1}, N_1(\nu), N_1(\Omega)\big) \leq C.
\end{equation}
Then, there exists $C'>0$ depending only on $(\mathcal U, \mathcal V)$, on the $C^1$ norm of $H$ on $\mathcal W$ and on $C$, 
such that 
$
   \|d\psi_t\|_{C^0(\mathcal W)}\leq C'.$
   \end{enumerate}
\end{proposition}
\begin{preuve}
For any $t\in [0,1]$, let $\Omega_t:= \Omega + td\nu.$ Then
for every $t\in [0,1]$, by hypothesis and (\ref{Caligula}), $S(\Omega_t, \mathcal W)\geq \frac{1}{2}S(\Omega, \mathcal W)$ which is positive since $\overline {\mathcal W}$ is compact, 
so that $\Omega_t$ is symplectic on $\overline{\mathcal W}$. 
We are looking for a 1-parameter family of diffeomorphisms $(\psi_t)_{t\in [0,1]}$ of $\mathcal W$ such that  $\forall t\in [0,1], \psi_t^*(\Omega_t ) = \Omega$ over $\psi_t (\mathcal U)$. Derivating in time, 
and assuming that $(X_t)_{t\in [0,1]}$ is a vector field that generates $(\psi_t)_t$, 
we obtain
$ \partial_t \Omega_t +d(i_{X_t}\Omega_t) = 0,$ or 
$d(\nu + i_{X_t}\Omega_t )= 0.$
We now inverse now the procedure. Let $(X_t)_{t\in [0,1]}$ be a family of vector fields on $\mathcal W$  such that
\begin{equation}\label{Dalida}
 \forall t\in [0,1], \ \forall x\in \mathcal W, \ i_{X_t(x)}\Omega_t(x) = -\nu(x).
 \end{equation}
Since $\Omega_t$ is non-degenerate, $X_t$ is uniquely defined, smooth and  by~(\ref{Neron}), 
\begin{equation}\label{ixe}
\forall t\in [0,1], \ \|X_t\|_{C^0(\mathcal W)}\leq  \frac{2}{S(\Omega, \mathcal W)}\|\nu\|_{C^0(\mathcal W)}.
 \end{equation} 
Let $\chi: \mathcal W\to [0,1]$ be a smooth cut-off function such that $\chi_{\overline {\mathcal V}}= 1$ and $\chi $ has support in $\mathcal W$. 
Let $(\psi_t)_{t\in [0,1]} $ be the 1-parameter family of diffeomorphisms associated to $\chi X_t$. 
By Lemma~\ref{majphivar},
\begin{equation}\label{bakounine}
\forall t\in [0,1], \ d(\psi_t, \id) \leq   \frac{2t}{S(\Omega, \mathcal W)}\|\nu\|_{C^0(\mathcal W)}.
\end{equation}
By hypothesis on $\|\nu\|$, this implies that $\forall t\in [0,1], \  \psi_t(\mathcal U)\subset \mathcal V.$
Since $\chi = 1$ over 
$\mathcal V$, we obtain 
$\psi_1^*\Omega_1 = \Omega$ over $\psi_1(\mathcal U)\subset \mathcal V$.

We now assume that~(\ref{stirner}) is satisfied
and want a bound for the derivative of $\psi_1$. Derivating equation (\ref{Dalida}) gives 
$\forall t\in [0,1], \ i_{X_t}\nabla \Omega_t + i_{\nabla X_t}\Omega_t = -\nabla \nu$
over $\mathcal W$, 
%$$i_{X_t}\nabla \Omega_t+ i_{\nabla X_t}\Omega_t = -\nabla %\nu,$$
so that 
$$ \max_{t\in[0,1]}\|\nabla X_t\|_{C^0(\mathcal W)}\leq \frac{2}{S(\Omega, \mathcal W)}(\| \nabla \Omega_t\|_{C^0(\mathcal W)} \|X_t\|_{C^0(\mathcal W)} + \|\nabla \nu\|_{C^0(\mathcal W)})$$
and
$ \max_{t\in[0,1]} \|\nabla (\chi X_t)\|\leq 
\max_{t\in [0,1]} \|\nabla X_t\|_{C^0}+\|d\chi\|\|X_t\|_{C^0}.$
By Lemma~\ref{majphivar} and~(\ref{bakounine}), this implies that
$\max_{t\in [0,1]}\|d\psi_t\|_{C^0} \leq C',$
where $C'$ depends only on the derivative of the metric on $\mathcal W$, on $C$ and $\chi$, hence on $(\mathcal V,\mathcal W)$. 
\end{preuve}

In Corollary~\ref{perturbation} below, we apply the latter proposition to the situation that is of interest for us: the construction of a symplectomorphism $\Psi$  between a symplectic submanifold $\Sigma$  in a ambient manifold $(M,\omega)$ and another submanifold $\phi(\Sigma)$ equipped the restriction of another symplectic structure $\omega +d\mu$ which is close to $\omega$. 
Then, the proof of Proposition~\ref{defo} will be a direct consequence of Corollary~\ref{perturbation}.
% % % % % % %
\begin{corollary}\label{perturbation} Under the hypotheses of Proposition~\ref{defo},
\begin{enumerate}
\item  there exists
a smooth isotopy  of embeddings $(\Psi_t)_{t\in [0,1]} : \Sigma\cap W \to \phi(\Sigma)$ satisfying
\begin{eqnarray*}
\forall t\in [0,1], \ \Psi_t^*\big((\omega+td\mu)_{|T\phi(\Sigma)}\big)&= &\omega_{|T\Sigma} \text{ on } U\cap \Sigma
\text{ with }
\Psi_t(U\cap \Sigma) \subset W\cap \phi(\Sigma),\\
\text{and }d(\Psi_1, \id_{|\Sigma\cap W}) &\leq & \frac{2}{S(\omega_{|T\Sigma}, W\cap \Sigma)}\| \phi^*(\lambda +\mu) -\lambda\|_{C^0(W)} +d(\phi, \id).
\end{eqnarray*}
\item If furthermore the hypotheses of~\ref{defo3}. of Proposition~\ref{defo} are satisfied, there exists $C'>0$ depending only on $U, V, \Sigma$, the $C^1$ norm of $h$ on $W$ and on $C$, such that
$\forall t\in [0,1],  \|d\Psi_t\|_{C^0(\Sigma\cap W)}\leq  C'.$
\end{enumerate}
\end{corollary}
% % % % % % % % % % %
\begin{preuve}
Let $j : \Sigma \to W $ be the natural injection, 
$ \Omega := j^* \omega
$ and $\nu :=j^* \big(\phi^*(\lambda+\mu)-\lambda\big),$
so that  $j^*\big( \phi^*(\omega+d\mu))= \Omega + d\nu.$
Choosing the metric $H$ on $\Sigma$ to be the induced one by the ambient metric $h$, the various estimates for $\nu$ are bounded
by the ones for $\phi^*(\lambda+\mu)-\lambda$, 
so that, using that the induced distance in $\Sigma$ is larger than 
the one in $M$, 
\begin{eqnarray*}
%N_1(\phi^*(\lambda+\mu)-\lambda, W)\leq C\\
\|\nu\|_{C^0(W\cap \Sigma )} \leq
\frac{1}{2}S(\omega_{|T\Sigma}, W\cap \Sigma)
\text{ dist}_\Sigma (U\cap \Sigma, \partial V\cap \Sigma) 
\text{ and } \|d\nu\|_{C^0(W\cap \Sigma)}\leq  
\frac{1}{2}S(\omega_{|T\Sigma}, W\cap \Sigma).
\end{eqnarray*}
By assertion~\ref{moser1} of Proposition~\ref{moser} applied to $(\Sigma, \Omega)$, $\mathcal W=\Sigma\cap W$,  $\mathcal V=\Sigma \cap V$, 
$\mathcal U=\Sigma \cap U$, 
$H=h_{|T\Sigma},$  and $\nu$, there exists  
a 1-parameter family of diffeomorphisms: $\psi_t: \Sigma\cap W \to \Sigma\cap W$ with support in $W\cap \Sigma$ such that
for any $t\in [0,1]$, $\psi_t(U\cap \Sigma)\subset V\cap \Sigma$, 
 \begin{equation}\label{zoroastre}
 d(\psi_t, \id)\leq \frac{2t}{S(\omega_{|T\Sigma}, W\cap \Sigma)}\|\nu\|_{C^0(W\cap \Sigma)}
 \leq \text{dist}(\mathcal U, \partial V),
 \end{equation}
 and
$\psi_t^*\Big(j^*\big(\phi^*(\omega+d\mu)\big)\Big)= \omega_{|T\Sigma} \text{ on } U\cap \Sigma.$
For any $t\in [0,1], $ let
$\Psi_t:= \phi\circ \psi_t: W\cap \Sigma \to Y.$ 
Then, since by hypothesis $d(\phi, \id)\leq \text{dist}(V, \partial W)$, we have by~(\ref{zoroastre}) $\Psi_t(U\cap \Sigma)\subset W\cap \phi(\Sigma)$. Moreover, 
 $$  \big(\Psi_t^*(\omega+d\mu)\big)_{|T\Sigma} = \Big(\psi_t^*\big(\phi^*(\omega+td\mu)\big)\Big)_{|T\Sigma} = \Omega \text 
 { on } U\cap \Sigma. $$
 This proves the first assertion. 

 Now, assume that the hypotheses of~\ref{defo3}. in Proposition~\ref{defo} are satisfied. Then
 $N_1(\nu, W\cap \Sigma) \leq  C$, so that by assertion~\ref{moser2} of Proposition~\ref{moser}, 
  there exists $C'>0$ depending only on $(\mathcal U, \mathcal V)$, the $C^1$ norm of $H$ and $C$, 
  $ \ \|d\psi_t\|_{C^0(W\cap \Sigma)}\leq C'$.
   Since $\|d\phi\|_{C^0(Y)}\leq C$, we have
  $\|d\Psi_t\|_{C^0(W\cap \Sigma)} \leq CC'$, 
hence the result after changing the definition of $C'$. 
\end{preuve} 
We can now give the proof of Proposition~\ref{defo}, which demonstrates the stability of a Lagrangian submanifold
in a symplectic submanifold when the latter and the symplectic form are perturbed.
\begin{preuve}[ of Proposition~\ref{defo}]
By Corollary~\ref{perturbation} there exists  a smooth diffeomorphism $$\Psi: \Sigma\cap U\to \Psi(\Sigma \cap U)\subset \phi(\Sigma)\cap W$$ such that $\big(\Psi^*(\omega+d\mu)\big)_{|T\Sigma} = \omega_{|T\Sigma}$ on $U\cap \Sigma$.
This implies that $\mathcal L' := \Psi(\mathcal L)$ is a smooth compact Lagrangian submanifold in $\big(\phi(\Sigma)\cap W, (\omega+d\mu)_{|T\phi(\Sigma)}\big)$. 

Assume now that the hypotheses of~\ref{defo2} in Proposition~\ref{defo} are satisfied. 
Then by Corollary~\ref{perturbation}, $d(\Psi, \id)\leq \frac{1}{4}\diam_M \mathcal L$. 
  Let $p,q\in \mathcal L$ such that $\diam_M (\mathcal L) = d_{M}(p,q).$
Then 
\begin{eqnarray*}
\diam_{\mathcal L'}(\mathcal L')\geq  d_{\mathcal L'} (\Psi(p), \Psi(q)) \geq d_{M} (\Psi(p), \Psi(q)) \geq
d_M(p,q) - 2d(\Psi, \id)
 \geq  \frac{1}{2}\diam_M(\mathcal L).
\end{eqnarray*}

Assume now that the hypothesis of~\ref{defo3} in Proposition~\ref{defo} is satisfied. Again by Corollary~\ref{perturbation},  
there exists $C'$ such that $\|d\Psi\|_{C^0(\Sigma\cap W)}\leq C'.$ This implies
$\diam_{\mathcal L'} (\mathcal L' )\leq C' \diam_{\mathcal L} (\mathcal L).$
Indeed, let $p',q'\in \mathcal L'$  and
$\gamma : [a,b]\to \mathcal L$ such that $\gamma$ is a shortest path in $\mathcal L$ between
$p:=\Psi^{-1}(p')\in \mathcal L$ and $q:=\Psi^{-1}(q')\in \mathcal L$. 
Then,
$$
d_{\mathcal L'} (p',q' )\leq 
Length_{\mathcal L'} (\Psi (\gamma))
=\int_a^b |d\Psi (\gamma) (\gamma'(t))| dt
 \leq  C'\diam (\mathcal L).
$$
\end{preuve}
\section{Proof of the main local theorem}\label{pot}

\subsection{The standard setting}

\begin{preuve}[ of Theorem~\ref{theorem2}]
  We adapt the barrier method of the real context in \cite{GWL}
  %~\cite[p. 1343]{NS1} 
  to our complex algebraic situation, and will use the quantitative Moser method given by Theorem~\ref{concrete}.
  For the reader's convenience, we begin by
the proof in the case of the standard random polynomials. 
Then we sketch the proof for the general 
setting of random holomorphic sections. 

Let $p\in \C[z_1, \cdots, Zz_n]$ regular such that $Z(p)\cap \mathbb B = \Sigma$. 
Since $p$ is regular, there exists $\eta>0$ such that $p : 2\mathbb B\to \C^r $ satisfies
the transversality condition 
\begin{equation}\label{tra}
 \forall z\in 2\mathbb B, \ |p(z)|< \eta\Rightarrow T\big(dp(z)\big)> \eta,
\end{equation}
where $T$ is defined by~(\ref{T}).

Since the probability measure is invariant under the symmetries of $\C P^n$, as well as the assertion of Theorem~\ref{theorem2}, it is enough to prove the theorem for 
$x=[1:0:\cdots : 0]$.
Let $z$ be the local holomorphic affine coordinates:
$$z=(z_1, \cdots, z_n):=  \Big(\frac{Z_1}{Z_0}, \cdots , \frac{Z_n}{Z_0}\Big)\in \C^n$$
defined on $\C P^n \setminus \{Z_0=0\}.$
Fix $\ep>0$ and let
$p_{\ep,d} (z) := p\Big(z\frac{\sqrt d}{\ep}\Big).$
Note that $Z(p_{\ep,d})= \frac{\ep}{\sqrt d}\Sigma$. 
Then for any $d\geq d(p)$,
let 
\begin{equation}\label{ped}
P_{\ep, d}(Z): =Z^d_0p_{\ep, d}\Big(\frac{Z_1}{Z_0}, \cdots , \frac{Z_n}{Z_0}\Big)\in \big(\C_{hom}^d [Z_0, \cdots , Z_n]\big)^r.
\end{equation}
  By construction,  $Z(P_{\ep, d})\subset \C P^n $ intersects the affine coordinate ball $B(0,\ep/\sqrt d)$ around $[1:0\cdots :0]$ along a small homothetical copy of $\Sigma $ and contains a copy of $\mathcal L$. Notice that $P_{\ep, d}$ is  singular, since $Z(P_{\ep, d})$  contains the hyperplane $\{X_0=0\}$ with multiplicity $d-d_0$.

In order to apply  the first item of Theorem~\ref{concrete}, 
we must have a bound for the perturbation of $\omega_0$ in $\omega$. 
For this, in our affine coordinates, let $\lambda_{FS} = d^c \log (1+|z|^2)$
and $\lambda_0 = d^c |z|^2$, 
that is 
$$ \lambda_{FS} = \frac{1}{2i}\frac{\sum_{i=1}^n z_i {d\bar z_i}-\bar z_i dz_i}{1+|z|^2}.$$
By definition $\omega_{FS} = d\lambda_{FS}$ and $\omega_0 = d\lambda_0$, 
so that 
$\lambda_{FS}= \lambda_0+ O(\|z\|^3)$, and $d\lambda_{FS}= \omega_0+O(\|z\|^2). $
Let $\psi $ the linear map $ \psi(z)=z\frac{\ep}{\sqrt d}$. Then, 
there exists a universal constant $K>0$  such that 
 the pull-backs $\lambda=\frac{d}{\ep^2} \psi^*\lambda_{FS}$ and  $\omega= \frac{d}{\ep^2}\psi^*\omega_{FS}$ satisfy
$$\forall d\geq d_{\mathcal L}:= \frac{16K}{\min(1, \diam_{Z(p)} \mathcal L)}, \  \|\lambda-\lambda_0\|_{C^0(2\mathbb B)} +\| \omega - \omega_0\|_{C^0(2\mathbb B)}\leq  K\frac{\ep^2}{d}\leq \frac{1}{16}\min (1, \diam_{Z(p)} \mathcal L),$$  
which is the bound needed in Theorem~\ref{concrete} for the perturbation form $\mu$ and its
differential, see conditions~(\ref{kikou}) and~(\ref{cloclo}).

Now let $H_P := P_{\ep,d}^\perp$ be the orthogonal space
to $P_{\ep, d}$ in $\big(\C_{hom}^d[Z], \langle \ , \rangle\big).$
We  use a decomposition
for our random polynomials adapted to $P_{\ep,d}$ an $H_P$.
Since the random polynomial can be written in any fixed orthormal basis, we can decompose our random polynomial $P$
 as 
 \begin{equation}\label{decomposition}
P=  a\frac{P_{\ep,d}}{\|P_{\ep, d}\|_{L^2}}+ R,
 \end{equation}
  where 
 $a$ is a complex Gaussian variable and $R\in H_P$ 
 is a Gaussian random polynomial for the induced law on $H_P$ and independent of $a$.
 The $L^2$-norm of $P_{\ep, d}$ is computed by Lemma~\ref{norm} below. 
 We want to prove that with uniform positive lower bound, $R$ does not perturb too much the first term, so that $P$
  still vanishes on a hypersurface diffeomorphic to $\Sigma$. 
 Hence, we need to know when  the vanishing locus of a perturbation of a function gives a diffeotopic perturbation of the vanishing locus of the function. 
For this, for any $d\geq d(p)$, we apply Theorem~\ref{concrete} to 
 \begin{eqnarray*}
 \forall z\in 2\mathbb B, \ f(z) := aP_{\ep, d}(1, z\frac{\ep}{\sqrt d})=ap(z) 
 \text{ and }g(z) :=\|P_{\ep, d}\|_{L^2} R(1, z\frac{\ep}{\sqrt d}).
  \end{eqnarray*}
By~(\ref{tra}) we have
 $$\forall a\in \C^*, \ \forall z\in 2\mathbb B, \ |f(x)|< |a| \eta\Rightarrow T\big(df(x)\big)>|a|  \eta.$$
We want now to give a uniform lower bound for the probability that the pair
of random functions $(f,g)$ on $2\mathbb B$ satisfies
the various conditions of Theorem~\ref{concrete}. 
In order to control the peturbation $g$, 
 we decompose it as $$g =p_1+p_2:= \frac{1}{2}(g+f) + \frac{1}{2}(g-f).$$
 Note that the law of $p_1:= g+f$ is the same of 
 $r(z) := \|P_{\ep, d}\|_{L^2} P(1,z \frac{\ep}{\sqrt d}),$ where $P$ follows the Fubini-Study law. The same holds for $p_2:=g-f$. 
 We use the trivial  inequality
 \begin{equation}\label{borne}
 \E \sup_{2\mathbb B} |g| \leq
 \frac{1}2\big( \E \sup_{2\mathbb B} |p_1|+\E  \sup_{2\mathbb B} |p_2|\big)\leq \E \sup_{2\mathbb B}|r|,
  \end{equation}
 and similarly for the average of the derivative of $g$. 
 Hence, it is enough to bound from above the norms
 of a random  $q$.

By Markov inequality, by the independence between $f$ and $p_1$, $p_2$, 
by the bound~(\ref{borne}), by Remark~\ref{norm2} and by Lemma~\ref{razorbak}, there exists $K_P>0$ depending only on $P$  that for
all  $0< \ep \leq 1$,  $d\geq d(p)$, $F>0$, $0<\alpha\leq 1$,  and  $\forall j\in \{0,1, 2\}$,
\begin{eqnarray}\label{chouchou}\nonumber
  \mathbb P_d\Big[ \exists a\in \C^*, \|g\|_{C^1}\leq |a|\eta/8, 
c_j(|a|\eta, f,g)\leq \alpha
\Big]&\geq  &\\ \nonumber
  \mathbb P_d\left[F\leq  |a|,  \|g\|_{C^1}\leq F\eta/8, 
\|g\|_{C^{j}}\leq F\alpha \frac{\eta^{2j+2}}{\|p\|^{2j+1}_{C^{3}}}\right]&\geq&  \\
\frac{1}{\pi}\int_{F< |a|} e^{-|a|^2}|da|
\Big( 1-\frac{K^2_P}{\ep^{2d(p)}\alpha^2 F^2}\Big).
\end{eqnarray}
Recall that $c_j$ is defined by~(\ref{cocorico}). For  $j=2$, let $ F=F_\epsilon:= 2\frac{K_P}{\alpha \ep^{d(p)}}$
and
$\alpha=\frac{1}{16}\diam_{\R^{2n}} \mathcal L$.
Then, there exists a constant $C_P>0$ depending only on $P$  such that for all $d\geq d(p)$, $0<\epsilon\leq 1$, 
the probability~(\ref{chouchou}) is bounded from above by 
$C_P\exp(-\frac{C_P}{\ep^{2d(P)}\diam^2_{\R^{2n}} \mathcal L}).$
By assertions \ref{con1}., \ref{con2}. and \ref{con3}. of Theorem~\ref{concrete}, 
there exists $C''$ depending only on $\diam_{\R^{2n}}\mathcal L$, such that for $d\geq \max \big(d(p), d_{\mathcal L}\big),$ 
with the same probability,
there exists  a topological  ball $B$ satisfying
$\frac{1}{2}\mathbb B\subset B\subset \frac{3}{2}\mathbb B$, and 
$\mathcal L'$ a compact smooth Lagrangian submanifold of $\big(Z(f+g)\cap B, \omega_{|Z(f+g)}\big)$,
satisfying 
$$ (\mathcal L, Z(f)\cap \mathbb B)  \sim_{diff} (\mathcal L', Z(f+g)\cap B)$$ 
with
$ \frac{1}{2}\diam_{Z(f)} (\mathcal L)\leq \diam_{\mathcal L'} (\mathcal L')
\leq \mathcal C'' \diam_{\mathcal L} (\mathcal L).$ Here, the metrics are the various restrictions of the standard metric $g_0$ on the ball. However, the push-forward of the metric $g_\omega$ on the unit ball by the coordinates $z\epsilon/\sqrt d$ converges uniformly in $0\leq \epsilon \leq 1$ to $g_0$ when $d$ grows to infinity. This implies the theorem.
\end{preuve}

\subsection{The general K\"ahler setting} 
The generalization of the proof of Theorem~\ref{theorem2} to random holomorphic sections
holds onto the concept of \emph{peak sections}, as in~\cite{GWL} and~\cite{GW}.  This object was used by Tian~\cite{tian} to give estimates for the Bergman kernel, and by S. K. Donaldson~\cite{donaldson96} to prove the existence of codimension 2 symplectic submanifolds. In a way, 
they were already used by H\"ormander to solve the Levi problem for Stein manifolds~\cite[Theorem 5.1.6]{hormander} .
They are used in the parallel paper~\cite{GD} for a deterministic proof of Corollary~\ref{coroX}. 

Let $(n, r, X,L,E,h_L,\omega, g_\omega, h_E, d\text{vol}, (\mathbb P_d)_{d\geq 1})$
be an ample probabilistic model
A peak section of $L^{\otimes d}$ at $x\in X$ is a holomorphic section which norm decreases
exponentially fast outside $x$, and almost vanishing at scale $\gg 1/\sqrt d$, 
like $X_0^d$ in the standard projective case near the point $[1:0 \cdots 0]$.
One of their crucial interest lies in the fact that a given peak section times the monomials~(\ref{monomials}) in
normal holomorphic coordinates form asymptotically an orthonormal family, 
which make the general
K\"ahler situation locally very similar to the standard projective one. 
\begin{preuve}[ of Theorem~\ref{theorem3}]
Let $x\in X$, and $e$ a local holomorphic trivialization of $L$ near 
$x$ such that $\|e\|_{h_L} $  is locally maximal at $x$, with $\|e(x)\|_{h_L} = 1$. 
Then there exists a uniform (in $x\in X$) constant $c>0$ such that  for any $y$ in a fixed neighborhood of $x$, 
\begin{equation}\label{pic}
\|e^{\otimes d} (y)\|_{h_L} \leq \exp(- cd\|x-y\|^2).
\end{equation}
This is implied by the fact that the curvature of $h_L$ is a K\"ahler form and the uniformity is implied by the compacity of $X$. 
Again, $X_0^d =e^d$ in the standard case. 
Let $(e_1, \cdots, e_r)$ be a local holomorphic trivialization of $E$, orthonormal at $x$. 
Then, $(e_1\otimes e^d, \cdots, e_r\otimes e^d)$ is a local holomorphic trivialization of $E\otimes L^d$ whose coordinates are \emph{called peak sections} for $x$. Now, let $(p_1, \cdots, p_r)$ be a polynomial map that defines the complex algebraic hypersurface $\Sigma$, and 
$$ s_{\ep, d, p}:= \Big( p_i \big(\cdot \frac{\sqrt d }{\epsilon} \big) e_i\Big)_{1\leq i\leq r} \otimes e^d,$$
which is a section of $E\otimes L^d$ defined in a fixed neighborhood of $x$, and is the equivalent of $P_{\ep, d}$, see~(\ref{ped}) in the standard case. 
Now by the H\"ormander $L^2$-estimates, see~\cite{hormander} or  \cite{skoda} for a bundle version of it,
$s_{\ep, d, p}$ can be perturbed in a global section $\sigma_{\ep, d, p}\in H^0(X,E\otimes L^d)$. Moreover, this is a classical result
in H\"ormander theory that 
the $C^1$ error produced by the perturbation on $B(x, \frac{\log d}{\sqrt d})$ is bounded by $\exp(-cd)$, see~\cite[Lemma 3.5]{GWL}. 
Here, the estimates~(\ref{pic}) are crucial. 
By Lemma~\ref{morse} this implies that $Z(s_{\ep, d, p})$ is 
a  complex $(n-r)$-submanifold which is an isotopic peturbation of $Z(p)$. 

The rest of the proof is very similar to the standard case. 
We decompose the random section $s\in H^0(X,E\otimes L^d)$ as
$$ s =a\frac{ s_{\ep, d, p}}{\|s_{\ep, d, p}\|_{L^2}}+ \rho,
$$
where $\rho \in s_{\ep, d, p}^\perp$ and $s_{\ep, d, p}^\perp$ 
is equipped with the restriction of the Gaussian measure, and $a$ follows a complex
normal law $N_{\C}(0,1)$. 
The $L^2$-norm of $s_{\ep, d, p}$ has a similar equivalent as $\|P_{\ep, d}\|_{L^2}$
given by Lemma~\ref{norm}.
Then we look the situation on $B(x,\ep/\sqrt d)$
which becomes a fixed $
\mathbb B\in \C^n $ after rescaling, and the sections are trivialized as functions with values in $\C^r$. 
Lemma~\ref{razorbak} still holds for the trivialization $q$ of the pertubation. In the proof of it, 
the only essential adaptation in the bundle case is the estimate~(\ref{RSA}),
where the modulus of the function is compared on $B(0, \ep/\sqrt d)$ with the Fubini-Study norm of it. 
In the present situation, a similar comparison holds, since the norm of $e^d$ varies only 
of a uniform positive multiplicative constant over $B(x,\ep/\sqrt d)$.

The Lagrangian part of the proof is the same, since the coordinates
at a point $x\in X$ we can choose holomorphic 
coordinates $z$ such that $\omega = z^*\omega_0$ at $x$, 
so that we can find a $1-$form $\lambda $
on the chart so that $\lambda - z^*\lambda_0 = O(|z|)$, 
which is the two only thing we need. 
\end{preuve}
%A peak section $\sigma_x$ for $x$ 
%is a local trivialization of $L^d$ 
%The following lemma establishes the existence of peak sections for $E\otimes L^d$. It is an almost direct
%consequence of~\cite[Lemma 1.2]{tian}:
%\begin{lemma}{\cite[Lemma 2.3.3]{GW}}\label{lemma tian}
%Let $(n,r, X, L, h_L, \omega, E, h_E) $ be a projective pack 
%$x\in X$, $(p_1, \cdots, p_n)\in \Nn^n$, $i\in\{1, \cdots, r\}$
%and $p'>0 $. There exists $d_0\in \Nn$ independent of $x$ such that for every $ d>d_0$,
%there exists $\si\in H^0(X,E\otimes L^d) $ with the property that $\|\si\|_{L^2(h_L)}= 1$ and
% if $ (z_1, \cdots, z_n)$ are  local holomorphic coordinates in the neighbourhood 
%of $x$ such that $\partial_{z_1}, \cdots \partial_{z_n}$ are orthonormal for $g_{\omega}$,  and $(e_1, \cdots e_k)$ is a local  holomorphic trivialization of $E$
%orthonormal at $x$,  we can assume that in a neighbourhood of $x$,
%\begin{equation}\label{sigma}
% \si(z_1, \cdots, z_n) = \lambda e_i\otimes e^d (1+ O(d^{-2p'} )) + O(\lambda |z|^{2p'}), 
% \end{equation}
%where 
%$ \lambda^{-2}= \int_{B(x,\frac{\log d }{\sqrt d})}
%\|e^d\|^2_{h_L^d} dV_{h_L}$, 
%with $dV_{h_L}= \omega^n/{\int_X \omega^n}$
% and where $e$ is a local trivialization of $L$ 
%whose potential $ -\log h_L (e,e)$ reaches a local minimum at $x$ with
%Hessian $\pi \omega(.,i.)$.
%\end{lemma} 
\begin{remark}
Instead of peak sections, we could have
use  the Bergman kernel, the Schwartz kernel for the projection onto the space of holomorphic sections, and the 2-point correlation function for our random model. This kernel  converges at scale $1/\sqrt d$ to a universel kernel, the \emph{Bargmann-Fock} kernel, see~\cite{BSZ}, which explains why the results on standard Fubini-Study random polynomials are similar to those for random holomorphic sections. This universality can be proved by peak sections, see~\cite{tian}. The kernel approach has the virtue that parts of
the proof can be adapted to other random models. However, we must not overestimate this interest for some reasons. Firstly, the fact that the zeros of the sections of given degree have the same  topology and symplectomophism type is very dependent on holomorphicity, or at least asymptotic holomorphy in the Donaldson~\cite{donaldson96} and Auroux settings~\cite{auroux}. Secondly, the projective hypersurface inherit a natural symplectic form, which is rarely the case for other models. Thirdly,
  the barrier method is very adapted to explicit local sections, like peak sections. Fourthly, the fact that this model is particularly well suited for polynomials is not directly seen from the kernel and need some asymptotic computation.  Lastly, let us notice that the Bergman kernel beween $x$ and $y$ is essentially represented by the value of the peak section associated to $x$ evaluated at $y$. 
\end{remark}
We finish this section with the proof that the smooth vanishing loci have all the same symplectomorphism type:
\begin{proposition}\label{nuke}
Let $1\leq r\leq n$ be an integer and $E\to X$ be a 
holomorphic vector bundle of rank $r$, $L\to X$ be a holomorphic line bundle equipped with a metric $h$ with positive curvature 
$-i\omega$. 
For any \emph{degree} $d\geq 1$, denote by $H_{reg}^0(X,E\otimes L^{ d})$ the space of holomorphic sections of $E\otimes L^{\otimes d}$which vanish transversally. 
Then for any $d$ large enough, 
$$ \forall (s,t) \in H_{reg}^0(X,E\otimes L^{ d})^2, 
\big(Z(s), \omega_{|Z(s)}\big)\sim_{symp} \big(Z(t), \omega_{|Z(t)}\big).$$
\end{proposition}
\begin{preuve}
First, by Bertini's theorem~\cite[p. 137]{griffiths}, $H^0_{sing}:= H_{reg}^0(X,E\otimes L^{ d})\setminus H_{reg}^0(X,E\otimes L^{ d})$ is of real codimension at least 2 in $H^0$. 
%To see this, fix a basis $\{S_0, \cdots, S_{N_d}\}$ of $H_0(X,E\otimes L^{ d})$ and define
%the incidence variety:
%$$ \mathcal I_d:= \{ (x, \lambda)\in X\times \C P^{N_d}, \sum_{i=1}^{N_d}\lambda_i S_i (x) = 0\}.$$
%Let $\pi : \mathcal I_d\to \C P^{N^d}$
%the natural projection on the second factor. Then, it is easy to see that $\pi$ is holomorphic and its set of critical values $\mathcal C\subset \C P^{N_d}$  correspond to the singular sections. Since by the Kodaira theorem, for $d$ large enough there are regular sections, $\mathcal C\not= \C P^{N_d}$. Since $\pi$ is holomorphic, $\mathcal C$ is a compact complex analytic subvariety of $\C P^{N_d}$, so that it has at least complex codimension 1. 
%
This implies that any pair $(s,t)\in H^0_{reg}$ are joined by a path of sections in $H^0_{reg}$. By Ehresmann theorem, this implies that $Z(s) $ is diffeomorphic to $Z(t)$. Now, for a continuous family of sections $(s_t)_{t\in [0,1]}$ in $H^0_{reg}$, 
since $\omega$ is a curvature of a line bundle, as its restriction to $Z(s_t)$, 
$[\omega_{|Z(s_t)}]\in H^2(Z(s_t), \mathbb Z)$. Consequently, the  pullback in $H^2(Z(s_0), \Z)$  of $[\omega_{|Z(s_t)}]$ by the diffeomorphism 
$ \psi_t : Z(s_0)\to Z(s_1)$ given by the former argument is constant. In other terms, $\psi_t^*[\omega_{|Z(s_t)}] =[\omega_{|Z(s_0)}]$. Then, the Moser argument (see~\cite[Theorem 3.17]{mcduff}) implies that the zero sets are in fact symplectomorphic. 
\end{preuve}

\subsection{Some simple lemmas}
In this paragraph we give the proofs of elementary and technical  lemmas
that are used in the core of the proof of the quantitative Moser deformation Proposition~\ref{defo}.

\paragraph{Lemmas for the deformations}
\begin{lemma}\label{majphi} Let  $m\geq 1$ be an integer, $(X_t)_{t\in [0,1]}$ be a $C^2$  family of 
vector field on $\R^n$ with compact support, and $(\phi_t)_{t\in [0,1]}$ be the associated flow. Then,
\begin{enumerate}
\item 
$\forall t\in [0,1], \ \|\phi_t - \id\|_{C^0(\R^m)} \leq 
 t\max_{t\in [0,1]} \|X_t\|_{C^0(\R^m)}  .$
\item Let $0\leq j\leq 2$ and $C>0$ be such that $\max_{t\in [0,1]}N_j(X_t, \R^m)\leq C$. Then, there exists $C'$ depending only on $C$ such  that
$$\forall t\in [0,1], \ \|\phi_t - \id\|_{C^j(\R^m)} \leq 
C' t\max_{t\in [0,1]} N_j(X_t)  .$$
\end{enumerate}
\end{lemma}
\begin{preuve}
First, it is classical that $\phi_t$ is $C^k$ in $(t,x)$.
We have 
$$\forall (x,t)\in \R^n\times [0,1], \ \phi_t(x) - x= \int_{0}^t X_s(\phi_s(x))ds,$$
which implies $\|\phi_t-\id\|_{C^0(M)}\leq \max_{t\in [0,1]} \|X_t\|_{C^0(M)}$
and 
$$d\phi_t -\id = \int_0^t d_x X_s(\phi_s(x)) \circ d\phi_sds.$$
Consequently, $\|d\phi_t-\id\|_{C^0}
\leq \max_t \|dX_t\|_{C^0} \big(t+
\int_0^t \|d\phi_s-\id \| ds\big).$
By Gronwall, this implies
\begin{equation}\label{grongron}
 \|d\phi_t-\id\|_{C^0}\leq t \max_t \|dX_t\|_{C^0} 
\exp\big(\max_t \|dX_t\|_{C^0}\big)\leq t e^C N_1(X) . 
\end{equation}
Now, 
$  d^{2} (\phi_t-Id) = d^{2}\phi_t =\int_0^t 
d^{2}_x X_s(\phi_s) d\phi_s\otimes d\phi_s +d_xX_s \circ d^{2}\phi_s ds.$
Together with estimate~(\ref{grongron}), this implies
$$ \|d^2\phi_t \|\leq \max_t \|d^{2}X_t\|_{C^0} (1+Ce^{C} )^2 +\max_t \|dX_t\|_{C^0}\int_0^t \|d^{2}\phi_s\|ds$$
so that by Gronwall,
$
\|d^{2}\phi_t\|
\leq  \max_t \|d^{2}X_t\|_{C^0} (1+Ce^{C} )^2 \exp(C).
$
\end{preuve}
% % % % % % % % % % %
Unfortunately, for manifolds we need a simplier version of the latter affine lemma. 
\begin{lemma}\label{majphivar} Let $(M,h)$ 
be a smooth Riemannian manifold,  
$(X_t)_{t\in [0,1]} $ be a $C^k$ family of  vector fields with compact support  $N$ and $(\phi_t)_{t\in [0,1]}$ the 1-parameter group of diffeomorphism generated by $(X_t)_t$. 
 Then
 \begin{enumerate}
\item 
$ \forall t\in [0,1], \ d(\phi_t, \id) \leq t\max_{t\in [0,1]} \|X_s\|_{C^0(M)}.$
\item 
Let $C>0$ be such that $\max_{t\in [0,1]}N_1(X_t, \R^m)\leq C$.
Then, there exists $C'$ depending only on $C$ and the $C^1$ norm of the metric on $N$, such that 
$\max_{t\in [0,1]} \|d\phi_t\|_{C^0(M)}
\leq C'.$
\end{enumerate}
\end{lemma}
\begin{preuve}
Again, it is classical that  $\phi_t$ is $C^k$ in $(t,x)$, and 
$$\forall (x,t)\in M\times [0,1], \  d(\phi_t(x), x) \leq Length ( \{ \phi_t(x)\}_{t\in[0,1]})\leq t\max_{s\in [0,1]} \|X_s\|_{C^0(M)}.$$
Let $x\in M$ in a local chart. If $t$ is small enough, \emph{In coordinates}, we have 
$$\phi_t(x) - x= \int_{0}^t X(\phi_s(x), s)ds,$$
so that 
$d\phi_t -\id = \int_0^t d_x X(\phi_s(x), s) \circ d\phi_sds.$
%$$d(\phi_t -\id) = \int_0^t d_x X(\phi_s(x), s) \circ (d\phi_s-\id)+d_xX(\phi_s(x), s) ds.$$
Then, there exists a constant $C$ depending
only on the compact support of $X$ and the $C^1$ norm of the metric $h$ in the coordinates, such that for any vector $Y\in \R^n$, 
$$|d\phi_t(Y)|_{\phi_t(x)} \leq C|Y|_{x} \big(1+\max_{s\in [0,t]} \|d_x X\|_{C^0(M)}\int_0^t\| d\phi_s\|_{\phi_s(x)}ds\big)$$
which implies
$\|d\phi_t\|_{\phi_t(x)} \leq C\big(1+\max_{s\in [0,t]} \|d_x X\|_{C^0(M)}\|\int_0^t\| d\phi_s\|_{\phi_s(x)}ds\big)$
and by Gronwall
$$ \|d\phi_t\|_{C^0}\leq C \exp(C\max_{[0,1]}\|d_xX\|),$$
so that there exists another constant
$C'$ depending on the chart, such that 
$$ \|d\phi_t\|_{C^0}\leq C \exp(C'\max_{[0,1]}N_1(X)\|).$$
Since we can cover the support of $X$ by a finite number of charts, this implies the result. 
%and
%$$ d^{2}\phi_t = d^{2} (\phi_t-Id) = \int_0^t 
%d^{2}_x X(\phi_s, s) d\phi\otimes d\phi +d_xX \circ d^{(2)}\phi_s ds.$$
%The latter implies
%$$ \|d^{(2)}\phi_t\|\leq \int_0^t \|d^{(2)}X\|_{C^0} \|d\phi_s\|_{C^0}^2 +\|dX\|\|d^{2}\phi_s\|ds .$$
%so that by Gronwall,
%\begin{eqnarray*}
%\|d^{(2)}\phi_t\|
%&\leq & t\|d^{(2)}X\|_{C^0} \|d\phi_t\|_{C^0}^2\exp \|dX\|\\
%&\leq &  t\|d^{(2)}X\|_{C^0}\big(1+\|d_xX\| \exp(\|d_xX\|)\big)^2
%\exp \|dX\|.
%\end{eqnarray*}
%%Assume that for $j\leq k-1$, 
%%$$d^{(j)}(\phi_t-\id) =\sum_{p=0, q\leq p}^j a_{p,j}\int_0^t d^{(p)}_x X_s(\phi_s(x)) \circ \otimes_{\ell=1}^p d^{\ell}(\phi_t-\id)
%%$$
\end{preuve}
The following Lemma~\ref{proj}  was used for  the proof 
 of last assertion~(\ref{phi}) of Proposition~\ref{morse}.
 \begin{lemma}\label{proj} Let $m\geq 1$, $1\leq p\leq m$ be two integers,  and 
 $\Phi: M(p,m)\times \R^m\to \R^m$, 
 where for any $(A,Y)\in M(p,m)\times \R^m$, $\Phi( A, Y) $ denotes the orthogonal projection of the origin onto 
 the space $\{X\in \R^m, AX= Y\}$. Then, for any $0\leq j\leq 2$ and 
 for any $(A,Y)\in M(p,m)\times \R^m$ such that $A$ is onto,
 \begin{eqnarray*}
  \Phi(A,Y) &\leq  & T^{-2}(A) \|A\||Y|,\\
  \|d_A \Phi(A,Y)\|&\leq& 3 T^{-4}(A)\|A\|^2 |Y| \text{ and }
  \|d_Y \Phi(A,Y)\|\leq T^{-2}(A) \|A\|,\\
  \|d_{A^2}^2 \Phi(A,Y)\|&\leq&  14\|A\|^3 T(A)^{-6}|Y| \text{ and }
\|d^2_{AY}\Phi \|\leq  3T(A)^{-4}\|A\|^2, 
  \end{eqnarray*}
  where $T $ has been defined by~(\ref{T}).
 \end{lemma}
 \begin{preuve}
 Write $A=(A_1, \cdots, A_p)^t$, where $A_k \in \R^m$ are column vectors. Since $A$ is onto, 
 $(\ker A)^\perp = \langle A_1, \cdots, A_p\rangle$
 and there exist a unique $\lambda(A,Y) = (\lambda_1, \cdots, \lambda_p)\in \R^p$, 
 such that $\Phi(A,Y) = \sum_{i=1}^p \lambda_i A_i \in A^{(-1)}(\{Y\})$, 
 which means  $AA^t \lambda = Y$, 
 so that 
 $$ \Phi(A,Y) = \big[K(A)Y, A\big],$$
 where $[\lambda, A]:= \sum_{i=1}^p \lambda_i A_i$
 and $K(A) :=  (AA^t)^{-1}$.
 This implies  that $\Phi$ is a smooth near $(A,Y)$
 for any $A$ onto, and linear in $Y$. 
 Since $$K(A)\leq T(A)^{-2},$$ 
 by Cauchy-Schwarz  
   $ |\Phi(A,Y)|\leq  T(A)^{-2}\|A\||Y|.$
 Derivating gives
  \begin{equation}\label{crochet}
   d\Phi(A,Y) = \big[dK(A)Y+KdY , A]+[KY, dA]
\end{equation}
so that 
$\|d_Y\Phi \|\leq 
T(A)^{-2}\|A\|.
$
If $Q(A):= AA^t$, then $\|dQ(A)\|\leq 2\|A\|
$ and 
$$\forall B\in M_{p,m}(\R), \ dK(A)(B) =- K dQ(A) B K,$$
so that, using $T(A)\leq \|A\|$, 
 \begin{eqnarray*}
 \|dK(A)\|&\leq & 2T(A)^{-4}\|A\|\\
\text{and }\|d_A\Phi\|& \leq & 2T(A)^{-4}\|A\|^2|Y| +T(A)^{-2}|Y|\leq 
3 T(A)^{-4}\|A\|^2|Y|.
\end{eqnarray*}
Now 
  \begin{equation}\label{crochet}
  d^2 \Phi(A,Y) = [d^2K(A) Y, A]+\text{Sym }\Big(
  [dK (A)Y, dA]+[K(A)dY,dA]
  +[dK(A)dY,A]\Big),
\end{equation}
where Sym means that the bracket is symmetrized in the two vectors in $M_{p,m}(\R)\times \R^m$ on which applies $d^2\Phi(A,Y)$.
Since $\ d^2K(A)= 
\text{Sym } KdQKdQK-Kd^2QK,$
we obtain
\begin{eqnarray*}
 \|d^2K(A)\|&\leq &8\|A\|^2 T(A)^{-6}+ 2T(A)^{-4} \leq
10\|A\|^2 T(A)^{-6}, \\
\|d_{A^2}^2 \Phi(A,Y)\|&\leq &
10\|A\|^3 T(A)^{-6}|Y|+4T(A)^{-4}\|A\||Y|\leq
14\|A\|^3 T(A)^{-6}|Y|\\
\text{and }
\|d^2_{AY}\Phi \|&\leq & 3T(A)^{-4}\|A\|^2.
\end{eqnarray*} 
 \end{preuve}

\paragraph{Lemmas for the barrier methods. }
\begin{lemma}\label{norm}
Let  $0< \ep\leq 1$, $p\in (\C[z_1, \cdots, z_n])^r$ and $P_{\ep, d}= Z_0^d p\big(\frac{\sqrt d}{\ep} \cdot\big)$.
Then, uniformly in $\ep$, 
$$ \|P_{\ep, d}\|_{L^2} \underset{d\to\infty}{\sim}  \frac{1}{\pi d^{n/2}}  \|p\big(\frac{\cdot}{\ep}\big)\|_{BF},$$
where $\|p\|^2_{BF}:= \frac{1}{\pi^{n}}  \int_{\C^{n}}  |p( y)|^2e^{-|y|^2}|dy|$
defines the Bargmann-Fock norm of $p$. 
\end{lemma}
\begin{remark}\label{norm2}
Note that for any $p$, 
there exists a constant $c>0$ such that for any $0<\ep\leq 1$, and $d\geq 1$, 
$ \|P_{\ep, d}\|_{L^2}\leq  \frac{c}{d^{\frac{n}{2}}\ep^{\deg p}}.$
\end{remark}
\begin{preuve}
We have, by definition of the Fubini-Study measure on $\C P^n$, 
$$ \|P_{\ep, d}\|_{L^2}^2= \int_{\mathbb S^{2n+1} }|P_{\ep, d}|^2 \frac{d\sigma}{2\pi},$$
where $d\sigma$ is the canonical measure on the sphere with volume 1 and the $2\pi$
factor corresponds to the volume of the fiber $U(1)$ of the quotient $S^{2n+1}\to \C P^n$. 
Since 
$$ \int_{\mathbb \C^{n+1} }|P_{\ep, d}|^2e^{-\|Z\|^2}  |dZ| = 
\int_{0}^\infty r^{2d+2n+1} e^{-r^2}dr
\int_{\mathbb S^{2n+1} }|P_{\ep, d}|^2 d\sigma = (d+n)! 
\int_{\mathbb S^{2n+1} }|P_{\ep, d}|^2 d\sigma,
$$
where $|dZ|$ denotes the Lebesgue measure on $\C^{n+1}$,
we have
$$ \|P_{\ep, d}\|_{L^2}^2 = \frac{1}{2\pi(d+n)!}\int_{\C^{n+1}} |Z_0|^{2d} \Big|p\big(\frac{\sqrt d}{\ep} \frac{Z'}{Z_0}\big)\Big|^2e^{-\|Z\|^2} |dZ|,$$
where $Z'= (Z_1, \cdots, Z_n)$. 
We use the change of variable $(W_0, w) = (Z_0,  \frac{Z'}{Z_0})$,  then $(w_0,w) = (W_0\sqrt{1+|w|^2}, w)$, and finally $y=\sqrt d w$ so that
\begin{eqnarray*}
 \|P_{\ep, d}\|_{L^2}^2 &=&  \frac{1}{2\pi(d+n)!\pi^{n+1}}\int_{\C^{n+1}} |W_0|^{2(d+n)} |p( \frac{\sqrt d}{\ep} w)|^2e^{-|W_0|^2(1+\|w\|^2)} |dW_0| |dw|\\
 & =&  \frac{1}{2\pi(d+n)!\pi^{n+1}} \int_{\C^{n+1}} |w_0|^{2(d+n)} e^{-|w_0|^2} |p( \frac{\sqrt d}{\ep} w)|^2\frac{1}{(1+\|w\|^2)^{d+n+1}} |dw_0||dw|\\
 & =& \frac{1}{d^{n}} \frac{1}{2\pi(d+n)!\pi^{n+1}} \int_{\C} |w_0|^{2(d+n)} e^{-|w_0|^2} |dw_0| \int_{\C^{n}} |p(\frac{y}\ep)|^2\frac{1}{(1+\frac{1}d\|y\|^2)^{d+n+1}} |dy|\\
 &\underset{d\to\infty}{\sim} & \frac{1}{d^{n}}  \frac{1}{\pi^{n+1}}
 \int_{\C^{n}}  |p(\frac{y}\ep)|^2e^{-|y|^2}|dy|
 \end{eqnarray*}
 uniformly in $\ep\leq 1$.
 
 Note that for any $(i_0, \cdots, i_n)\in \Nn^{n+1}$ such that $\sum_k i_k = d$, 
 \begin{eqnarray*}
   \|Z_0^{i_0}\cdots Z_n^{i^n}\|^2_{L^2}&=& \frac{1}{2\pi (d+n)!}  \int_{\mathbb \C^{n+1} }\Pi_{k=0}^{n+1}|Z_k^{i_k}|^2e^{-\|Z\|^2}  |dZ| \\
   & =& \frac{1}{2\pi (d+n)!}  \Pi_{k=0}^{n+1} \int_{\mathbb \C }|z|^{2i_k}e^{-|z|^2}  |dz|\\
   & = & \frac{1}{(d+n)!}  \Pi_{k=0}^{n+1} \int_{0}^\infty r^{2i_k+1}e^{-r^2} dr = \frac{i_0! \cdots i_n!}{(d+n)!}
 \end{eqnarray*} 
\end{preuve}
The next lemma was proved in a real and general K\"ahler situation in~\cite{GWL}. We give a proof of it in the polynomial setting, in order the article to be self-contained and simple. 
 \begin{lemma}\label{razorbak} Let $1\leq r\leq n$ be integers, $\ep>0$, 
$R\in (H_{d,n+1})^r$ be a random polynomial mapping of maximal degree $d$ and $q(z)= R(1,z\frac{\ep}{\sqrt d})$, where $z= (z_1, \cdots, z_n)$. Then, there exists $C>0$ depending only on $n$ and $r$ such that 
for any $d\geq 1$, any $0<\epsilon\leq 1$, any $0\leq j\leq 2$, 
  $\mathbb E \big( \sup_{2\mathbb B} |d^jq|^2\big)\leq  C_n\frac{(d+n)!}{d!} $.
 \end{lemma}
 \begin{preuve}
 Since $q$ is holomorphic, we can 
use the mean value inequality for plurisubharmonic functions applied to $|q|^2$ (see~\cite{hormander}):
 $$ \forall z\in 2\mathbb B,\  |q(z)|^2 \leq \frac{1}{\text{Vol } \mathbb B}\int_{z+\Bb} |q|^2(u) du,$$
 so that 
$  \mathbb E \big(\sup_{2\mathbb B} |q|^2\big) 
\leq 
\frac{1}{\text{Vol } \mathbb B}\int_{3\Bb} \mathbb E (|q(u)|^2)du.
$
We have by (\ref{poly})
\begin{equation}\label{RSA}
\forall z\in 2\mathbb B, \  \mathbb E (|q(z)|^2)
=
\mathbb E \big(|R(1, z\frac{\ep}{\sqrt d})|^2\big) = 
\mathbb E \big(\|R\|_{FS}^2(1, z\frac{\ep}{\sqrt d})\big)\big(1+\frac{|z|^2\ep^2}{d}\big)^{2d}. 
\end{equation}
By definition of the measure, $\mathbb E \big(\|R\|_{FS}^2)$
is constant over $\C P^n$. Remembering that the coordinates of $R$ 
are independent random polynomials, 
we obtain (see decomposition~\ref{poly}),
$$\mathbb E \big(\|R\|_{FS}^2)(1, z\frac{\ep}{\sqrt d})= \mathbb E \big(\|R\|_{FS}^2(0)) = r
\mathbb E \Big(\frac{(d+n)!}{d!} |a_0|^2\Big)= r\frac{(d+n)!}{d!}. $$
Moreover
$ \forall d\geq 1, \ \forall z\in 2\mathbb B, \ \big(1+\frac{|z|^2\ep}{d}\big)^{2d}\leq e^{18\ep^2},$
hence the first estimate of the Lemma. 

For the second estimate, define for any holomorphic function $f= (f_1, \cdots, f_r) : \C^n \to \C^r$
 $$\|df\|_2^2 := \sum_{i=1}^r \sum_{j=1}^n \Big|\frac{\partial f_i}{\partial z_j}\Big|^2.$$
 Notice that $\|df\|\leq \|df\|_2$, where $\|df \|$ is the operator norm used in Proposition~\ref{morse}.
We have, similarly to the first estimate since the complex derivatives of $q$ are holomorphic, 
$$  \mathbb E \big(\sup_{2\mathbb B} \|dq\|_2^2\big) 
\leq 
\frac{1}{\text{Vol } \mathbb B}\int_{3\Bb} \mathbb E \big(\|dq(u)\|^2_2\big)du.
$$
with $\|dq(u)\|_2^2 = \frac{\ep^2}{d }\|d_{Z'} R(1, u\frac{\ep}{\sqrt d})\|^2_2,$
where $Z'=(Z_1, \cdots, Z_n)$.
As before 
$$ \mathbb E \big(\|d_{Z'} P(1, u\frac{\ep}{\sqrt d})\|^2_2\big) \leq
\mathbb E \big(\|d_{Z'} P\|_{FS}^2(1, 0)\big)e^{18\ep^2}
$$
with, using the linear part in $Z'$ of the decomposition~(\ref{poly}), 
$$  \mathbb E \big(\|d_{Z'} P(1, 0)\|_{FS}^2\big) = 
r\sum_{i=1}^n \mathbb E \Big(\frac{(d+n)!}{(d-1)!1!} |a_{0\cdots 1\cdots 0}|^2\Big) = rn\frac{(d+n)!}{(d-1)!}$$
which implies the second estimate of the Lemma. The last estimate is similar.
\end{preuve}

\bibliographystyle{amsplain}
\bibliography{Lagrangian.bib}

\noindent Univ. Grenoble Alpes, Institut Fourier \\
F-38000 Grenoble, France \\
CNRS UMR 5208  \\
CNRS, IF, F-38000 Grenoble, France
\end{document}